\documentclass[10pt]{amsart}

\usepackage{amssymb} 
\usepackage[mathscr]{eucal} 
\usepackage{xypic}
   \CompileMatrices 
\usepackage{amsmath} 
\usepackage{epsfig}

\newtheorem{TEO}{Theorem}[section]
\newtheorem{PROP}[TEO]{Proposition}

\newtheorem{DEF}[TEO]{Definition}
\newtheorem{COR}[TEO]{Corollary}
\newtheorem{REM}[TEO]{Remark}

\newcommand\dual{\mathrel{\raise3pt\hbox{$\underline{\mathrm{\thinspace d
\thinspace}}$}}}

\newcommand\proj{\mathbb P}
\newcommand\Z{\mathbb Z}

\newcommand\R{\mathbb R}

\newcommand\Co{\mathbb C}

\def\Z{{\mathbb Z}}
\def\R{{\mathbb R}}
\def\Z{{\mathbb Z}}

\title{Real Kodaira Surfaces}
\author{Paola Frediani}
\thanks{ Research performed in the realm of the  
{\em EU Project "E.A.G.E.R." }
}
\date{}

\begin{document}

\begin{abstract}
In this paper we give the topological classification of real 
primary Kodaira surfaces and we describe in detail the structure of the 
corresponding moduli space. 

Moreover, we use the notion of the
orbifold fundamental group of a real variety, which was also the main tool in the classification of real hyperelliptic surfaces achieved in \cite{cf}. 

Our first result  is  that if $(S, \sigma)$ is a real
primary Kodaira surface, then the differentiable type of the pair $(S,
\sigma)$ is completely determined by the orbifold fundamental group
exact sequence. \\
This result allows us to determine all the possible topological types of
$(S, \sigma)$.

Finally,  we show that once we fix the topological type of $(S, \sigma)$
corresponding to a real primary Kodaira surface, the corresponding
moduli space is irreducible (and connected).

\vspace{0.2cm}

 
\vspace{0.2cm}

\end{abstract}

\maketitle

\section*{Introduction}

The purpose of this paper is to achieve the topological classification of real Kodaira surfaces and to describe the structure of their moduli space. 

The main tool that we used in the topological and differentiable classification is the notion the orbifold fundamental group exact sequence of a real variety, whose relevance in real geometry has been pointed out in \cite{cf}. 

To be more precise, a smooth real variety is a pair  $(X, \sigma)$, consisting of the data of a
smooth complex manifold $X$ of complex dimension
$n$ and of an antiholomorphic involution $\sigma: X
\rightarrow X$ (an involution  $\sigma$ is a map whose square is the
identity).

$X$ is a complex manifold, so it is determined by a differentiable manifold $M$ and a complex structure $J$ on the complexification of the real tangent bundle of $M$.

If we consider the same manifold $M$ together with the complex
structure $-J$, we obtain a complex manifold  which is called the conjugate
of $X$ and denoted by $\bar{X}$.

The involution $\sigma$ is now said to be antiholomorphic if it provides an
isomorphism between the complex manifolds $X$ and $\bar{X}$ (and then $(X,
\sigma)$ and
$(\bar{X}, \sigma)$ are also isomorphic as pairs).

If $(X, \sigma)$ is a compact real variety, one would like to describe the isomorphism classes of the pairs $(X, \sigma)$, or the possible topological or differentiable types of the pairs $(X, \sigma)$.

We notice that already the problem of describing the topological type of the real part of $X$, $X({\R}) := Fix(\sigma)$ can be rather difficult.

Recall that Hilbert's 16-th problem is a special case of the last
problem but for the more general case of a pair of real varieties
$(Z \subset X , \sigma)$.

For a smooth real variety, we have the quotient double covering
$\pi : X \rightarrow Y = X/<\sigma>$, and the quotient $Y$ is called the Klein
variety of
$(X, \sigma)$.

In dimension $n=1$ the datum of the Klein variety is equivalent to the datum
of the pair $(X, \sigma)$, but this is no longer true in higher dimension,
where we will need also to specify the covering $\pi$.

The covering  $\pi$ is  ramified on the real part of
$X$, namely,
$X': = X({\R}) = Fix(\sigma)$, which
   is either empty, or a real submanifold of real dimension $n$.

If $X': = X({\R}) = Fix(\sigma)$ is empty,
   the orbifold fundamental group of
$Y$ is just defined as
   the fundamental group of $Y$.

If $ X' \neq \emptyset $, we may take a fixed point $x_0 \in
Fix(\sigma) $ and observe that $\sigma$ acts on the fundamental group
$\pi_1(X, x_0)$: we can therefore define the orbifold fundamental group as
the semidirect product of the normal subgroup $\pi_1(X, x_0)$ with the cyclic
subgroup of order two generated by $\sigma$. It is easy to verify then that
changing the base point does not alter  the isomorphism class of the
following exact sequence, yielding the orbifold fundamental group as an
extension

$$ 1 \rightarrow \pi_1(X) \rightarrow \pi_1^{orb}(Y) \rightarrow {\Z}/2
\rightarrow 1 $$

(changing  the base point only affects the choice of a splitting of the above
sequence).

In \cite{cf} we studied real hyperelliptic surfaces and we proved that the topological and differentiable type of a real hyperelliptic surface is completely determined by the orbifold fundamental group exact sequence. Furthermore, once we fix the topological type of a real hyperelliptic surface, the corresponding moduli space is irreducible and connected. 

We also claimed that the orbifold fundamental group exact sequence is a powerful
topological invariant of the pair
$(X, \sigma)$ in the case where $X$ has large fundamental group, in particular in the case where $X$ is a $K(\pi,1)$.

Primary Kodaira surfaces are non K\"ahler surfaces which are $K(\pi, 1)$ and Kodaira gave a description of them as quotients of ${\Co}^2$ by a group acting by affine transformations (cf. \cite{koi}). So we thought that we could try to study real primary Kodaira surfaces, in order to give another issue (besides real hyperelliptic surfaces) where the topology of the pair
$(S, \sigma)$ is determined by the orbifold fundamental group exact sequence.

Our first result is  

\begin{TEO} Let $(S, \sigma)$ be a real primary Kodaira surface. Then the
differentiable type of the pair $(S, \sigma)$ is completely determined by the
orbifold fundamental group exact sequence.
\end{TEO}

Primary Kodaira surfaces are non K\"ahler surfaces of Kodaira dimension $0$, their first Betti number is three and their first homology group is isomorphic to ${\Z} \oplus {\Z} \oplus {\Z} \oplus {\Z}/m{\Z}$, where $m \in {\Z}$, $m \geq 1$.

Moreover the torsion coefficient $m$ completely determines the differentiable type of a Kodaira surface.

The second result on the topology of real Kodaira surfaces that we have is the following.

\begin{TEO}
Let us fix the topological type for a Kodaira surface $S$, i.e. we fix the torsion coefficient $m \in {\Z}$ of the first homology group of $S$. If $m \equiv 0 \ (mod \ 2)$, then the number of topologically different real Kodaira surfaces is equal to 17; if $m \equiv 1 \ (mod \ 2)$, then the number of topologically different real Kodaira surfaces is equal to 13.
\end{TEO}

As a consequence we obtain

\begin{COR}

Let $(S, \sigma)$ be a real Kodaira surface. Then the real part
   $S(\mathbb R)$ is either
  the empty set, or  
   a disjoint union of $t$ tori, where $ 1 \leq t \leq 4$.  
\end{COR}

Regarding then the complete description of the moduli space of real Kodaira surfaces, we have the following main result that asserts that the differentiable type of the pair $(S, \sigma)$ 
determines the deformation type. This is true for real K\"ahler surfaces of Kodaira dimension less or equal to 0, but it is  false  already for complex surfaces
if the  Kodaira dimension  equals 2,  cf. \cite{ca1},\cite{ca2},
\cite{ma}, \cite{ca5}, \cite{kha-ku}).

\begin{TEO} Fix the topological type of $(S, \sigma)$ corresponding to a real
Kodaira surface. Then the moduli space of the real surfaces $(S',
\sigma ')$ with the given topological type is irreducible (and connected).
\end{TEO}

To prove the last result we used the description of the moduli space of complex Kodaira surfaces of a given differentiable type given by Borcea in \cite{bo1}.

\begin{TEO}
(\cite{bo1})
The moduli space corresponding to isomorphism
classes of complex structures on a fixed topological (differentiable)
type $S_0$ of a Kodaira surface may be identified to the product of
the complex plane with a punctured disk.
\end{TEO}

Concerning now the Enriques classification of real algebraic surfaces, it has
been focused up to now mostly on the classification of the topology of the
real parts, the topological classification of real rational surfaces going
back to Comessatti (\cite{co1}\cite{co2}\cite{co3}), as well as the
classification of real abelian varieties (\cite{co3}, see also \cite{si},
\cite{s-s}). A complete description of the deformation classes of real structures on minimal ruled surfaces is given by Welschinger (\cite{we}).

In the case of real $K3$ - surfaces we have the classification by Nikulin and
Kharlamov (\cite{ni},
\cite{kha}), for the real Enriques surfaces the one by Degtyarev and Kharlamov
(\cite{dekha1}, \cite{dekha3}).

For real hyperelliptic surfaces we have already mentioned the paper \cite{cf}.

Finally, partial results on real ruled and elliptic surfaces have been
obtained by Silhol (\cite{si}) and by Mangolte (\cite{man3}).

The  paper is organized as follows:

in section $1$ we recall the description given by Kodaira of (primary) Kodaira surfaces and the results of Borcea on the moduli space of Kodaira surfaces of a given topological type.

In section $2$ we describe the possible liftings of an antiholomorphic involution $\sigma$ on $S$ to the universal covering ${\Co}^2$ of $S$. In particular it turns out that all such liftings can be represented as affine transformations of ${\R}^4$.

In section $3$ we first recall the notion of the orbifold fundamental group of a real variety and we show that the representation of the
orbifold fundamental group as a group of affine transformations of ${\R}^4$
is uniquely determined, up to isomorphism, by the abstract structure of the
group.

Then, we show that, once this affine representation is fixed, the moduli
space for the compatible complex structures is irreducible and connected.

Finally we describe all the possible topological types of a real Kodaira surface. 

In section $4$ we explain how to determine the topology of the real part of a real Kodaira surface and we give a list of all the possible real parts as disjoint unions of tori.

\vspace{1cm}

{\bf Acknowledgements.} I would like to thank Prof. F. Catanese for the several interesting conversations on the subject and in particular on the possibility of using the orbifold fundamental group in order to topologically classify real Kodaira surfaces.

I also would like to thank Prof. L. Badescu for having pointed out the paper \cite{bo1} to my attention during a very nice stay at the IMAR in Bucharest in March 2002.

\vspace{0.2cm}

\section{Basics on Kodaira surfaces}

Kodaira proved the following theorem (cf. \cite{koi})

\begin{TEO}
Let $S$ be compact complex smooth surface. If the canonical bundle of
$S$ is trivial, then $S$ is a $K3$ surface, a complex torus or an
elliptic surface of the form ${\Co}^2 /G$, where ${\Co}^2$ is the
space of the two complex variables $(z_1, z_2)$ and $G$ is a properly
discontinuous group of affine transformations without fixed points of
${\Co}^2$ which leave invariant the $2$-form $dz_1 \wedge dz_2$. The
first homology group of the elliptic surface ${\Co}^2/G$ is 
$$H_1({\Co}^2/G, {\Z}) \cong {\Z} \oplus {\Z} \oplus {\Z} \oplus {\Z}/m.$$
\end{TEO}

The surfaces $S = {\Co}^2/G$ in the theorem are called primary Kodaira
surfaces and they admit a holomorphic locally trivial fibration over
an elliptic curve with an elliptic curve as typical fibre. 

In \cite{koi} Kodaira also proved that the fundamental group $G$ of $S$
can be generated by the elements $g_1, g_2, g_3, g_4$, which as
covering transformations have the form 
$$g_j (z_1,z_2) = (z_1 +\alpha_j, z_2 + \bar{\alpha_j} z_1 + \beta_j),$$
with 
$$\alpha_1 = \alpha_2 =0,$$
$$\bar{\alpha_3} \alpha_4 - \bar{\alpha_4} \alpha_3 = m \beta_2, \ m \in {\Z}, \ m \neq 0, $$
$\beta_1, \beta_2$ linearly indipendent over ${\R}$, i.e. $\beta_1 \bar{\beta_2} - \bar{\beta_1} \beta_2 \neq 0$.
We have $g_j = A_j z + b_j$, where 

$$A_j = 
\left( \begin{array}{cc} 
1 & 0 \\ 
\bar{\alpha_j} & 1
\end{array} \right)
$$
and 
$$
b_j = \left( \begin{array}{c} 
\alpha_j\\ 
\beta_j
\end{array} \right)
$$

\begin{REM}
The centre $Z$ of $G$ is the subgroup generated by $g_1$ and $g_2$, $Z = <g_1, g_2> \cong {\Z}^2$.

The commutator $[g_3, g_4] = g_3g_4 g_3^{-1} g_4^{-1} = g_2^{m}$.

Therefore we have a central extension
$$1 \rightarrow Z \cong {\Z}^2 \rightarrow G \rightarrow G/Z \cong
{\Z}^2 \rightarrow 1 \ \ \ (*)$$

\end{REM}

{\bf Proof.}

Let $g(z) = Az +b$, $\gamma(z) = \Lambda z + \delta$, then $g^{-1} (z) = A^{-1} z - A^{-1} b$.
We observe that $A_j \in 
\left\{
\left( \begin{array}{cc} 
1 & 0 \\ 
* & 1
\end{array} \right)
\right\} \cong ({\Co}, +)$
therefore 
$$[g,\gamma](z) = z + (Id - \Lambda)b - (Id -A)\delta.$$

Now the proof is an easy computation.
\hfill Q.E.D. \\

Set $E_{\beta} = {\Co}/({\Z}\beta_1 + {\Z} \beta_2)$ and $\pi: {\Co}^2 \rightarrow {\Co}^2/Z= {\Co} \times E_{\beta}$ be the projection. 
We have then 
$${\Co}^2/Z = {\Co} \times E_{\beta} \rightarrow {\Co}^2/G =({\Co} \times E_{\beta})/{\Z}^2,$$
and the map

$${\Co}^2/G =({\Co} \times E_{\beta})/{\Z}^2 \rightarrow {\Co}/{\Z} \alpha_3 + {\Z} \alpha_4$$
$$[(z_1, z_2)] \mapsto [z_1]$$
has fibre $E_{\beta}$.

One can prove that from a differentiable
viewpoint, Kodaira surfaces are completely determined by the torsion
coefficient $m$ of their first integral homology group (cf. \cite{bo1}).

If $S$ is a compact complex surface with trivial canonical bundle, a non zero global holomorphic two-form $\eta$ satisfies 
\begin{equation}
\label{eta}
d \eta =0, \ \eta \wedge \eta =0, \ \eta \wedge \bar{\eta} >0 \ at \ every \ point \ of \ S.
\end{equation}

Conversely, let $S_0$ be the underlying differentiable manifold, then any global complex valued two-form $\eta$ satisfying (\ref{eta}) defines a complex structure on $S_0$ with respect to which $\eta$ is holomorphic and nowhere null. 

Let $p \in \proj (H^2(S_0, {\Co}))$ be a point corresponding to the cohomology class of a global holomorphic non zero two form for some complex structure on $S_0$ with a trivial canonical bundle, the $p$ lies in the open set $D$ determined on the quadric $p \cdot p =0$ by the condition $p \cdot \bar{p} >0$. Here the product is cup product on $H^2(S_0, {\Co})$. The group of orientation preserving diffeomorphisms on $S_0$ acts on $D$, and we have the following results of Borcea (\cite{bo1}).

\begin{TEO}
Given a line $p \in D$ of $H^2(S_0, {\Co})$ there exist representatives $\eta$ of $p$ satisfying conditions (\ref{eta}). 
\end{TEO}   

\begin{TEO}
Any two such representatives define isomorphic complex analytic structures on $S_0$.  

Furthermore any complex structure on $S_0$ occurs in this manner i.e. it has a trivial canonical bundle.
\end{TEO} 

Therefore he proves that a parameter space for isomorphism classes of complex structures on $S_0$ is the quotient of $D$ by the action of orientation preserving diffeomorphisms on $S_0$ and we have the following theorem.

\begin{TEO}
(\cite{bo1})
The moduli space corresponding to isomorphism
classes of complex structures on a fixed topological (differentiable)
type $S_0$ of a Kodaira surface may be identified to the product of
the complex plane with a punctured disk.
\end{TEO} 

Borcea shows that if we take coordinates $(x_1, y_1, x_2, y_2)$ on the universal covering of $S_0$, a basis of $H^2(S_0,{\Co})$ is given by the forms $\theta_{ij}$, $1\leq i<j\leq 4$, $(i,j) \neq (1,2), (3,4)$, where $\theta_{ij} = \omega_i \wedge \omega_j$, with 
$$\omega_1 = dx_1, \ \omega_2 = dy_1, \ \omega_3 = dx_2 - x_1 dx_1 - y_1 dy_1, \ \omega_4 = dy_2 - x_1 dy_1 + y_1 dx_1.$$

So a two form $\eta$ can be written as follows:
$\eta = p_{13} \theta_{13} +  p_{23} \theta_{23} + p_{14} \theta_{14} + 
p_{24} \theta_{24}$, where $(p_{13}, p_{23}, p_{14}, p_{24})$ are homogeneous coordinates on $\proj H^2(S_0, {\Co})$ and $D$ is given by:
$$p_{13} p_{24} - p_{23} p_{14}=0,$$
$$-p_{13} \bar{p}_{24} + p_{23} \bar{p}_{14} +p_{14} \bar{p}_{23} - p_{24} \bar{p}_{13} >0.$$
One can also show (see \cite{bo1}) that the action of the orientation preserving diffeomorphisms of $S_0$ on $D$ is the following 
\begin{equation}
\label{*}
 \left( \begin{array}{c} 
p_{13} \\ 
p_{23} \\
p_{14} \\
p_{24}
\end{array} \right)
\mapsto 
\left( \begin{array}{cc} 
M & -\frac{2k}{m} M \\
0 & e \cdot M 
\end{array} \right)
\left( \begin{array}{c} 
p_{13} \\ 
p_{23} \\
p_{14} \\
p_{24}
\end{array} \right)
 \end{equation}
where $M = \left( \begin{array}{cc} 
a & b \\ 
c & d
\end{array} \right)$, $a,b,c,d \in {\Z}$ $ad -bc = e = \pm 1$, $k \in {\Z}$, $m$ is the torsion coefficient of the first homology group $H_1(S_0, {\Z})$.
So one can easily see (cf. \cite{bo1}) that the map
$$(p_{13}, p_{23}, p_{14}, p_{24}) \mapsto (p_{14}/p_{24}, p_{13}/p_{14})$$
defines an isomorphism of $D$ onto $(H_+ \times H_+) \cup (H_- \times H_-)$, where $H_+$ and $H_-$ denote respectively the upper and lower half planes.
Finally $p_{14}/p_{24}$ undergoes a modular transformation, while $p_{13}/p_{14}$ is translated by $-2k/m$. The quotient of $H_+$ by the modular group is isomorphic to the complex plane, the quotient of $H_+$ by an infinite cyclic group of real translations is isomorphic to the punctured disk. 

Now if we set $a = e = -d = -1$, $b = c = k =0$, then $p_{14}/p_{24} \mapsto - p_{14}/p_{24}$, $p_{13}/p_{14} \mapsto - p_{13}/p_{14}$ and the two components are interchanged. \\

\section{Symmetries of Kodaira surfaces}

Let us now assume that $S$ is a real Kodaira surface, and let $\sigma: S \rightarrow S$ be an antiholomorphic involution. Then we can find a lifting $\tilde{\sigma}$ of $\sigma$ to the universal cover ${\Co}^2$.

\begin{PROP}
Let $\tilde{\sigma}$ be a lifting of $\sigma$ to ${\Co}^2$, then $ \tilde{\sigma}$ is an affine transformation. 
\end{PROP}
{\bf Proof.}
Since $\tilde{\sigma}$ is a lifting of $\sigma$ to the universal covering, $\tilde{\sigma}$ acts by conjugation on $G$, therefore it also acts on the centre $Z$ of $G$, because it is characteristic. 
 
Hence $\tilde{\sigma}$ induces an antiholomorphic map $\bar{\sigma} : {\Co}^2/Z = {\Co} \times E_{\beta} \rightarrow {\Co}^2/Z = {\Co} \times E_{\beta}$ as it is illustrated in the following diagram:

$$\diagram {\Co}^2 \dto  \rto^{\tilde{\sigma}}      & {\Co}^2 \dto  \\
            {\Co}^2/Z = {\Co} \times E_{\beta} \dto  \rto^{\bar{\sigma}} & {\Co}^2/Z = {\Co} \times E_{\beta} \dto\\
{\Co}^2/G = S  \rto^{\sigma} & {\Co}^2/G = S \\
\enddiagram$$

Therefore, for all $\gamma \in Z$ there exists a $\gamma' \in Z$ such that 
$\tilde{\sigma} \circ \gamma = \gamma' \circ \tilde{\sigma}$, 
so if we write $\tilde{\sigma}(z_1, z_2) = (\sigma_1(z_1,z_2), \sigma_2(z_1, z_2))$, we have 
$$\sigma_1(z_1, \gamma(z_2)) = \sigma_1(z_1, z_2),$$
$$\sigma_2(z_1, \gamma(z_2)) = \gamma'(\sigma_2(z_1, z_2)).$$
Since ${\Co}_{z_2}/Z \cong E_{\beta}$ which is compact, we immediately
see that $\sigma_1$ is constant in $z_2$, $\sigma_1(z_1, z_2) =
\sigma_1(z_1)$. Furthermore, since $Z$ acts by translations on $
{\Co}_{z_2}$, we see that $\sigma_2$ is affine antiholomorphic in
$z_2$ and we can write
\begin{equation}
\label{2}
\sigma_2(z_1, z_2) = a_2(z_1) \bar{z_2} + c_2(z_1).
\end{equation}
Now, if we write $\gamma(z_2) = z_2 + \delta$, with $\delta \neq 0$, $\gamma'(z_2) = z_2 +
\delta'$, from (\ref{2}) we obtain 
$\sigma_2(z_1, \gamma(z_2)) = a_2(z_1)(\bar{z_2} + \bar{\delta}) +
c_2(z_1) =  \gamma'(\sigma_2(z_1, z_2))= a_2(z_1)(\bar{z_2}) +
c_2(z_1) + \delta'$, so $a_2(z_1) \bar{\delta} = \delta'$ and
$a_2(z_1) = a_2 \in {\Co}$ is a constant.

$G/Z =:H \cong {\Z}^2$ acts by translations on
${\Co}_{z_1}$ and ${\Co}_{z_1}/H = E_{\alpha}$ which is
compact. $\tilde{\sigma}$ normalizes $H$, therefore for every $g \in
H$ there exists a $g' \in H$ such that 
$$\sigma_1(g(z_1)) = g'(\sigma_1(z_1)).$$
So we see that $\sigma_1$ is affine antiholomorphic in $z_1$ and we
can write
\begin{equation}
\label{1}
\sigma_1(z_1) = c \bar{z_1} + d
\end{equation}

Now we use the fact that there exists $g  \in G$ such that $\tilde{\sigma} \circ g_3  = g \circ \tilde{\sigma}$, and $g' \in G$ such that $\tilde{\sigma} \circ g_4  = g' \circ \tilde{\sigma}$.

We observe that the action of $G$ on the first component is given by translations, therefore we have

$$
\tilde{\sigma} \circ g_3(z_1, z_2) = (c(\bar{z_1} + \bar{\alpha_3}) +
d, a_2(\alpha_3 \bar{z_1} + \bar{z_2} + \bar{\beta_3}) + c_2(z_1 +
\alpha_3))= $$
$$
g \circ \tilde{\sigma} (z_1, z_2) = g ( c \bar{z_1} + d,
a_2 \bar{z_2} + c_2(z_1)) = 
\left( \begin{array}{cc} 
1 & 0 \\ 
x & 1
\end{array} \right)
\left( \begin{array}{c} 
c \bar{z_1} + d \\ 
a_2 \bar{z_2} + c_2(z_1)
\end{array} \right)
+ \left( \begin{array}{c} 
\bar{x}\\ 
*
\end{array} \right)
$$
Analogously for $g_4$ we have 

$$
\tilde{\sigma} \circ g_4(z_1, z_2) = (c(\bar{z_1} + \bar{\alpha_4}) +
d, a_2(\alpha_4 \bar{z_1} + \bar{z_2} + \bar{\beta_4}) + c_2(z_1 +
\alpha_4))= $$
$$
g' \circ \tilde{\sigma} (z_1, z_2) = g' ( c \bar{z_1} + d,
a_2 \bar{z_2} + c_2(z_1)) = 
\left( \begin{array}{cc} 
1 & 0 \\ 
x' & 1
\end{array} \right)
\left( \begin{array}{c} 
c \bar{z_1} + d \\ 
a_2 \bar{z_2} + c_2(z_1)
\end{array} \right)
+ \left( \begin{array}{c} 
\bar{x'}\\ 
*
\end{array} \right)
$$
So from the first components we find $x = \bar{c} \alpha_3$, $x' =
\bar{c} \alpha_4$. The second components yield
$$a_2 \alpha_3 \bar{z_1} + a_2 \bar{\beta_3} + c_2(z_1 + \alpha_3) = xc
\bar{z_1} + xd + c_2(z_1) + \ constants,$$ 
$$a_2 \alpha_4 \bar{z_1} + a_2 \bar{\beta_4} + c_2(z_1 + \alpha_4) = x'c
\bar{z_1} + x'd + c_2(z_1) + \ constants.$$ 
By derivation w.r.t. $\bar{z_1}$ we obtain 

$$a_2 \alpha_3  + \frac{\partial c_2}{\partial \bar{z_1}}(z_1
+ \alpha_3) = xc + \frac{\partial c_2}{\partial \bar{z_1}}(z_1),$$

$$a_2 \alpha_4  + \frac{\partial c_2}{\partial \bar{z_1}}(z_1
+ \alpha_4) = x'c + \frac{\partial c_2}{\partial \bar{z_1}}(z_1).$$
So $\frac{\partial c_2}{\partial \bar{z_1}}(z_1) = h \bar{z_1} + f$
and by subtituting this expression in the last two equations we get 
$$a_2 = \frac{xc - h \bar{\alpha_3}}{\alpha_3} = \frac{x'c - h
\bar{\alpha_4}}{\alpha_4}.$$
Since $\tilde{\sigma}^2 \in G$, we immediately see that $c \bar{c} =
1$.  Hence, since $ x = \bar{c} \alpha_3$, $x' = \bar{c} \alpha_4$, we have $xc = \alpha_3$, $x'c = \alpha_4$ and $a_2 = 1 - h
\frac{\bar{\alpha_3}}{\alpha_3} = 1 - h
\frac{\bar{\alpha_4}}{\alpha_4}$, so $h(\bar{\alpha_3} \alpha_4 -
\alpha_3 \bar{\alpha_4}) = 0$, that yields $h =0$, $a_2 =1$. 

Therefore we have $\frac{\partial c_2}{\partial \bar{z_1}}(z_1) = f$,
hence $c_2(z_1) = f \bar{z_1} + g$, and 
$$\tilde{\sigma}
\left( \begin{array}{c} 
z_1 \\ 
z_2
\end{array} \right)
= \left( \begin{array}{cc} 
c & 0 \\ 
f & 1
\end{array} \right)
\left( \begin{array}{c} 
\bar{z_1}\\ 
\bar{z_2}
\end{array} \right)
   + \left( \begin{array}{c} 
d\\ 
g
\end{array} \right).
$$

Let us now impose the condition ${\tilde{\sigma}}^2 \in G$. 
$$\tilde{\sigma}^2
\left( \begin{array}{c} 
z_1 \\ 
z_2
\end{array} \right)
= 
\left( \begin{array}{cc} 
|c|^2 & 0 \\ 
f\bar{c} + \bar{f} & 1
\end{array} \right)
\left( \begin{array}{c} 
z_1 \\ 
z_2
\end{array} \right)
+ 
\left( \begin{array}{c} 
c \bar{d} + d  \\ 
f \bar{d} + \bar{g} + g
\end{array} \right) =$$
$$
\left( \begin{array}{cc} 
1 & 0 \\ 
j\bar{\alpha_3} + s \bar{\alpha_4}  & 1
\end{array} \right)
\left( \begin{array}{c} 
z_1 \\ 
z_2
\end{array} \right)
+ 
\left( \begin{array}{c} 
j \alpha_3 + s \alpha_4  \\ 
*
\end{array} \right),
$$
where $j, s \in {\Z}$. 

Therefore, as we have already noticed we have $|c|^2 = 1$ and $c
\bar{f} + f = j \alpha_3 + s \alpha_4 = c \bar{d} + d$.
\hfill{Q.E.D.}

\begin{REM}
Let $\tilde{\sigma}$ be a lifting of ${\sigma}$ to the universal covering
${\Co}^2$ as above, then we have 
$$\tilde{\sigma}
\left( \begin{array}{c} 
z_1 \\ 
z_2
\end{array} \right)
= \left( \begin{array}{cc} 
c & 0 \\ 
f & 1
\end{array} \right)
\left( \begin{array}{c} 
\bar{z_1}\\ 
\bar{z_2}
\end{array} \right)
   + \left( \begin{array}{c} 
d\\ 
g
\end{array} \right)
$$
where $|c|^2 = 1$, $c\bar{f} + f = c \bar{d} + d \in \Gamma_{\alpha} =
\alpha_3 {\Z} + \alpha_4 {\Z}$ and $\bar{\Gamma}_{\beta} =
\Gamma_{\beta} = \beta_1 {\Z} + \beta_2 {\Z}$, $c
\bar{\Gamma}_{\alpha} = \Gamma_{\alpha}$. 

\end{REM}
{\bf Proof.}
It only remains to show that $\bar{\Gamma}_{\beta} =
\Gamma_{\beta}$ and $c \bar{\Gamma}_{\alpha} = \Gamma_{\alpha}$. 

The condition $\tilde{\sigma} g_i \tilde{\sigma}^{-1} \in Z$, $i =1,2$
gives $\bar{\Gamma}_{\beta} = \Gamma_{\beta}$, while by imposing
$\tilde{\sigma} g_j \tilde{\sigma}^{-1} \in G$, $j = 3,4$ we obtain $c
\bar{\Gamma}_{\alpha} = \Gamma_{\alpha}$.
\hfill{Q.E.D.}

\section{Topological types of real Kodaira surfaces}

First of all we recall the notion of the orbifold fundamental group exact sequence, that we have introduced in \cite{cf}.

Let $(X, \sigma)$ be a smooth real variety of  dimension
$n$ (i.e., $X$ is a smooth complex manifold of complex dimension $n$ given
together with an antiholomorphic involution $\sigma$). Then we have a 
double covering $\pi : X
\rightarrow Y = X/<\sigma>$ ramified on $X' = X({\R}) = Fix(\sigma)$. 
Set $Y' := \pi(X')$.

We will define the orbifold fundamental group exact sequence of $(X, 
\sigma)$ as the isomorphism class of a given extension

\begin{equation}
\label{orb}
1 \rightarrow {\pi}_1(X,x_0) \rightarrow {\pi}_1^{orb}(Y,y_0) 
\rightarrow {\Z}/2 \rightarrow 1
\end{equation}

The choice of a base point will however create some technical difficulties.

\begin{DEF}
\begin{enumerate}
\item
If $X' = \emptyset$, then we define $\pi_1^{orb}(Y,y_0) = \pi_1(Y,y_0)$.

\item If $X' \neq \emptyset$ and $x_0 \in X'$, $y_0 = \pi(x_0)$, 
$\sigma$ acts on $\pi_1(X,x_0)$ and we define $\pi_1^{orb}(Y,y_0)$ to 
be the semidirect product of the normal subgroup $\pi_1(X, x_0)$ with 
the cyclic
group of order $2$ generated by an element which will be denoted by 
${\tilde{\sigma}}_0$ and
whose action
  on $\pi_1(X, x_0)$ by conjugation is the one of $\sigma$.

\item If $n \geq 3$, $Y' \neq \emptyset$ and $x_0 \not\in X'$, define 
the orbifold fundamental group of $Y$ based on $y_0$ as $\pi_1(Y-Y', 
y_0)$.

\item Assume $x_0 \not \in X'$, $X' \neq \emptyset$, assume moreover 
$dim_{{\Co}} X = 2$. Then the orbifold fundamental
group of
$Y$ with base point $y_0 = \pi(x_0)$ is defined to be the quotient of 
$\pi_1(Y - Y', y_0)$ by the subgroup normally generated by 
${\gamma_1}^2, ..., {\gamma_m}^2$, where $Y'_1,...,Y'_m$ are the 
connected components of $Y'$ and $\gamma_i$ is a simple loop around 
$Y'_i$.
\end{enumerate}

Since in all cases we have a well defined exact sequence
$$
1 \rightarrow {\pi}_1(X,x_0) \rightarrow {\pi}_1^{orb}(Y,y_0) 
\rightarrow {\Z}/2 \rightarrow 1,
$$
this will be called the orbifold fundamental group exact sequence.
\end{DEF}

\begin{PROP}

The isomorphism class of the fundamental group exact sequence is 
independent of the choice of
$y_0$.
\end{PROP}
{\bf Proof.}
This is well known when comparing cases $(1)$, $(3)$ and $(4)$ which 
are mutually exclusive.

In case $(2)$ ($X' \neq \emptyset$) we claim that 
${\pi}_1^{orb}(Y,y_0)$ is independent of the choice
of $x_0 \in X'$.
In fact, let $\delta$ be a path connecting $x_0$ with $x_1$: then the map
$$\gamma \mapsto \delta^{-1} \gamma \delta$$
yields an isomorphism between $\pi_1(X,x_0)$ and $\pi_1(X,x_1)$.
The action of $\sigma$ on $\pi_1(X,x_1)$ reads out on $\pi_1(X,x_0)$ 
as the composition
$$\gamma \mapsto \delta^{-1} \gamma \delta \mapsto 
\sigma(\delta)^{-1} \sigma(\gamma) \sigma(\delta) \mapsto \delta 
\sigma(\delta)^{-1} \sigma(\gamma) \sigma(\delta) \delta^{-1}.$$

But this action is precisely the conjugation by ${\tilde{\sigma}}_1:= 
\delta \sigma(\delta)^{-1} {\tilde{\sigma}}_0$.

Since $ {\tilde{\sigma}}_1$ is an element of order $2$, we obtain 
that the split extensions
$$1 \rightarrow {\pi}_1(X,x_0) \rightarrow {\pi}_1^{orb}(Y,y_0) 
\rightarrow {\Z}/2 \rightarrow 1$$
and

$$1 \rightarrow {\pi}_1(X,x_1) \rightarrow {\pi}_1^{orb}(Y,y_1) 
\rightarrow {\Z}/2 \rightarrow 1$$
  are isomorphic.

To relate case $(2)$ with the other two it suffices, once $x_0 \in 
X'$ and $x_1$ are given, to choose a splitting of the extension

$$1 \rightarrow {\pi}_1(X,x_1) \rightarrow {\pi}_1^{orb}(Y,y_1) 
\rightarrow {\Z}/2 \rightarrow 1$$

simply by taking $\gamma_j$ with $j$ such that $y_0 \in Y'_j$.

\hfill{Q.E.D.}

\begin{REM}
Let us explain the definition of $\pi_1^{orb}(Y, y_0)$ in the case 
where $n \geq 2$ and $x_0 \not\in X'$.

$X'$ is a real submanifold of
codimension $n$, hence the map
$$ \pi_1(X - X', x_0) \rightarrow \pi_1(X,x_0)$$
is surjective if $n \geq 2$ and it is
an isomorphism for $n \geq 3$. The singularities of $Y$ are contained in $Y'
= \pi (X')$ and there we have a local model ${\R}^n \times ({\R}^n/(-1))$.
Therefore
$Y$ is smooth for $n = 2$ and topologically singular for $n \geq 3$.  The
local punctured fundamental group $\pi_1(Y - Y')_{loc}$ is isomorphic to
${\Z}$ for $n = 2$, while it is isomorphic to ${\Z}/2$ for $n \geq 3$. This
means that the kernel of the surjection
$$\pi_1(Y - Y',y_0) \rightarrow \pi_1(Y,y_0)$$ is normally generated by loops
$\gamma$ around the components of $Y'$. If
$n \geq 3$, then we automatically have $\gamma^2 = 1$.
\end{REM}

Let $\tilde{X}$ be the universal covering of $X$, so that $X =
\tilde{X}/\pi_1(X)$. The exact sequence (\ref{orb}) defines a group 
which is the
group of liftings of the action of ${\Z}/2 \cong \{Id_X, \sigma\}$ to 
$\tilde{X}$, so that $Y =
\tilde{X}/\pi_1^{orb}(Y)$. \\

\begin{REM}
If $X' \neq \emptyset$, then the exact sequence
$$1 \rightarrow \pi_1(X) \rightarrow \pi_1^{orb}(Y) \rightarrow 
{\Z}/2 \rightarrow 1$$
always splits, as follows by the definition.
\end{REM}

\begin{TEO}
\label{orfun}
Let $(S, \sigma)$ be a real Kodaira surface. Then the differentiable
type of $(S, \sigma)$ is completely determined by the orbifold
fundamental group exact sequence
\begin{equation}
\label{OFG}
1 \rightarrow G = \pi_1(S) \rightarrow \hat{G}  \rightarrow {\Z}/2 \rightarrow 1 
\end{equation}
\end{TEO}
{\bf Proof.}
First of all we see that the fundamental group exact sequence
determines the topological type of the surface $S$. 
In fact consider the homotopy exact sequence of the fibration $S
\rightarrow E_{\alpha}$ with fibre $E_{\beta}$
\begin{equation}
\label{G} 
1 \rightarrow Z \cong {\Z}^2 = <\gamma_1, \gamma_2> \rightarrow G \rightarrow G/Z \cong
{\Z}^2 = <\gamma_3, \gamma_4> \rightarrow 1 
\end{equation}

Kodaira proved that there exists a representation of $G$ in $A(2, 
{\Co})$ such that $\gamma_i$, $i = 1,2$ correspond to the translations
$g_i$, $i = 1,2$. Furthermore if $\pi: S \rightarrow E_{\alpha}$ is the
fibration, there exist $f_3, f_4 \in \pi_1(S)$ such that $\pi_*(f_3) =
\gamma_3$, $\pi_*(f_4) = \gamma_4$, and such that $f_3, f_4$ are
represented by the affine transformations $g_3$ and $g_4$.

By possibly taking other generators of the fundamental group of the
fibre, $\pi_1(E_{\beta}) =Z$, we can assume that we have the relation $g_3
g_4 g_3^{-1} g_4^{-1} = g_2^{m}$ and obviously (\ref{G}) determines the
integer $m$. 

Now, since the torsion coefficient $m$ determines
the differentiable type of $S$ (cf. \cite{bo1}), $\pi_1(S)$ determines $S$ differentiably.

Now we assume that $(S, \sigma)$ is real. Then we have seen that we
have a representation of a lifting $\tilde{\sigma}$ of $\sigma$
to the universal covering ${\R}^4$ of $S$ as an affine map of the form
 
$$\tilde{\sigma}  
\left( \begin{array}{c} 
x_1\\
y_1\\
x_2\\
y_2 
\end{array} \right)
= \left( \begin{array}{cccc} 
c_1 & c_2 & 0 & 0\\
c_2 & -c_1 & 0 & 0\\
f_1 & f_2 & 1 & 0\\
f_2 & -f_1 & 0 & -1 
\end{array} \right)
\left( \begin{array}{c} 
x_1\\
y_1\\
x_2\\
y_2 
\end{array} \right)+ 
\left( \begin{array}{c} 
d_1\\
d_2\\
\gamma_1\\
\gamma_2
\end{array} \right)
$$

The orbifold fundamental group exact sequence 
$$1 \rightarrow G \rightarrow \hat{G} \rightarrow {\Z}/2 \rightarrow 1$$
determines the orbifold fundamental group exact sequence of the action of $\sigma_1$ on
the elliptic curve $E_{\alpha}$. In fact we have already observed that
an element $\tilde{\sigma}$ of $\hat{G} - G$
acts by conjugation on $Z$, therefore it acts on $G/Z = H$. So we have determined an extension $1 \rightarrow H \rightarrow \hat{H} \rightarrow {\Z}/2 \rightarrow 1$, which is the orbifold fundamental group exact sequence for the real elliptic curve $E_{\alpha}$. Since for a real elliptic curve the orbifold fundamental group exact sequence determines the topological type (cf. e.g. \cite{cf}), we have shown that we can fix the topological type of
$(E_{\alpha}, \sigma_1)$, where $\sigma_1$ denotes as above the first
component of $\tilde{\sigma}$. 
We have thus three different topological cases for the action of
$\sigma_1$ on the universal covering ${\R}^2$ of $E_{\alpha}$, which can be distinguished as follows. 

We can assume that the linear part of $\sigma_1$, which is
given by the matrix 
$$\left( \begin{array}{cc} 
c_1 & c_2\\
c_2 & -c_1
\end{array} \right)$$ is one of the following 

A) $\left( \begin{array}{cc} 
c_1 & c_2\\
c_2 & -c_1
\end{array} \right) = \left( \begin{array}{cc} 
1 & 0\\
0 & -1
\end{array} \right)$ 

B) $\left( \begin{array}{cc} 
c_1 & c_2\\
c_2 & -c_1
\end{array} \right) = \left( \begin{array}{cc} 
0 & 1\\
1 & 0
\end{array} \right)$ 

If $Fix(\sigma_1) \neq \emptyset$ (which always happens in case B)),
we can assume $d_1=d_2 =0$. Otherwise if in case A) we have $Fix(\sigma_1) = \emptyset$, one can easily prove that we can assume $d_1 = (1/2)$ the $+1$-eigenvector of the linear part of $\sigma_1$, $d_2 =0$ (see for instance \cite{cf} lemma 5.5).

Every element $h \in G$ has the following form:

$$h  
\left( \begin{array}{c} 
x_1\\
y_1\\
x_2\\
y_2 
\end{array} \right)
= \left( \begin{array}{cccc} 
1 & 0 & 0 & 0\\
0 & 1 & 0 & 0\\
a & b & 1 & 0\\
-b & a & 0 & 1 
\end{array} \right)
\left( \begin{array}{c} 
x_1\\
y_1\\
x_2\\
y_2 
\end{array} \right)+ 
\left( \begin{array}{c} 
a\\
b\\
z\\
w
\end{array} \right)
$$

We would like now to determine the $f_i$'s and
the $\gamma_j$'s in the expression of $\tilde{\sigma}$. 
In order to do this we impose in cases A) and B) the condition that
for all $h \in G - Z$, $\tilde{\sigma} h \tilde{\sigma}^{-1} \in G-Z$.

Let us treat at first case A).

The linear part of $ \tilde{\sigma} h \tilde{\sigma}^{-1}$ is 
$$
\left( \begin{array}{cccc} 
1 & 0 & 0 & 0\\
0 & -1 & 0 & 0\\
f_1 & f_2 & 1 & 0\\
f_2 & -f_1 & 0 & -1 
\end{array} \right)
\left( \begin{array}{cccc} 
1 & 0 & 0 & 0\\
0 & 1 & 0 & 0\\
a & b & 1 & 0\\
-b & a & 0 & 1 
\end{array} \right) 
\left( \begin{array}{cccc} 
1 & 0 & 0 & 0\\
0 & -1 & 0 & 0\\
-f_1 & f_2 & 1 & 0\\
f_2 & f_1 & 0 & -1 
\end{array} \right)
=
$$
$$\left( \begin{array}{cccc} 
1 & 0 & 0 & 0\\
0 & 1 & 0 & 0\\
a & -b & 1 & 0\\
b & a & 0 & 1 
\end{array} \right)$$

Now we compute the translations parts and we obtain

$$
\left( \begin{array}{cccc} 
1 & 0 & 0 & 0\\
0 & 1 & 0 & 0\\
a & -b & 1 & 0\\
b & a & 0 & 1 
\end{array} \right) 
\left( \begin{array}{c} 
-d_1\\
0\\
-\gamma_1\\
-\gamma_2 
\end{array} \right)
+ $$
$$\left( \begin{array}{cccc} 
1 & 0 & 0 & 0\\
0 & -1 & 0 & 0\\
f_1 & f_2 & 1 & 0\\
f_2 & -f_1 & 0 & -1 
\end{array} \right)
\left( \begin{array}{c} 
a\\
b\\
z\\
w 
\end{array} \right)
+
\left( \begin{array}{c} 
d_1\\
0\\
\gamma_1\\
\gamma_2 
\end{array} \right)=
\left( \begin{array}{c} 
a\\
-b\\
f_2 b + f_1 a + z - a d_1\\
-b_1 d_1 + f_2 a -w -f_1 b
\end{array} \right). 
$$

Now, since 
 $\tilde{\sigma}h  \tilde{\sigma}^{-1} \in G$, we know the translation part of $\tilde{\sigma}h  \tilde{\sigma}^{-1}$, which must be of the form $$\left( \begin{array}{c} 
a\\
-b\\
\delta\\
\epsilon
\end{array} \right)
$$
so we know $\delta$ and $\epsilon$ and from the computation above we must have

$$\left( \begin{array}{cc} 
a & b\\
-b & a
\end{array} \right)
\left( \begin{array}{c} 
f_1\\
f_2
\end{array} \right)
=
\left( \begin{array}{c} 
\delta -z + a d_1\\
\epsilon +w +b d_1
\end{array} \right).
$$
We are assuming that $h \in G- Z$, so we know that $a^2 + b^2
\neq 0$, therefore $f_1$ and $f_2$ are uniquely determined, since the orbifold exact sequence determines the conjugation
action of $\tilde{\sigma}$ on $G$. 

Now we treat case B).

Here we may assume $d_1 =d_2 =0$. 
$$\tilde{\sigma}
\left( \begin{array}{c} 
x_1\\
y_1\\
x_2\\
y_2
\end{array} \right)= \left( \begin{array}{cccc} 
0 & 1 & 0 & 0\\
1 & 0 & 0 & 0\\
f_1 & f_2 & 1 & 0\\
f_2 & -f_1 & 0 & -1 
\end{array} \right)
\left( \begin{array}{c} 
x_1\\
y_1\\
x_2\\
y_2
\end{array} \right)+
\left( \begin{array}{c} 
0\\
0\\
\gamma_1\\
\gamma_2
\end{array} \right).
$$ 
Again, if $h \in G-Z$, we know that $\tilde{\sigma} h
\tilde{\sigma}^{-1} \in G -Z$, and therefore we know the action of 
$\tilde{\sigma} h \tilde{\sigma}^{-1}$ on the universal covering.  The linear part of
$\tilde{\sigma} h \tilde{\sigma}^{-1}$ is 
$$
\left( \begin{array}{cccc} 
0 & 1 & 0 & 0\\
1 & 0 & 0 & 0\\
f_1 & f_2 & 1 & 0\\
f_2 & -f_1 & 0 & -1 
\end{array} \right)
\left( \begin{array}{cccc} 
1 & 0 & 0 & 0\\
0 & 1 & 0 & 0\\
a & b & 1 & 0\\
-b & a & 0 & 1 
\end{array} \right)
\left( \begin{array}{cccc} 
0 & 1 & 0 & 0\\
1 & 0 & 0 & 0\\
-f_2 & -f_1 & 1 & 0\\
-f_1 & f_2 & 0 & -1 
\end{array} \right) = $$
$$
\left( \begin{array}{cccc} 
1 & 0 & 0 & 0\\
0 & 1 & 0 & 0\\
b & a & 1 & 0\\
-a & b & 0 & 1 
\end{array} \right)
$$ 
The translation parts are
$$ 
\left( \begin{array}{c} 
b \\
a \\
a f_1 + b f_2 +z\\
a f_2  - b f_1 -w 
\end{array} \right)=
\left( \begin{array}{c} 
b \\
a \\
x\\
y
\end{array} \right),
$$ where $\left( \begin{array}{c} 
b \\
a \\
x\\
y
\end{array} \right)$ is the translation part of $\tilde{\sigma} h \tilde{\sigma}^{-1}$.

So we have 
$$
\left( \begin{array}{cc} 
a & b\\
-b & a 
\end{array} \right)
\left( \begin{array}{c} 
f_1\\
f_2 
\end{array} \right)
= \left( \begin{array}{c} 
x -z\\
y+w 
\end{array} \right)
$$
and we are able to determine $f_1$ and $f_2$, since we know $h$ and
$\tilde{\sigma} h \tilde{\sigma}^{-1}$ and $a^2 + b^2 \neq 0$.  

So we can determine the linear part of $\tilde{\sigma}$ in all
cases. Now we want to find the translation part 
$ \left( \begin{array}{c} 
d_1 \\
d_2 \\
\gamma_1 \\
\gamma_2 
\end{array} \right)$ of which we already know the $d_i$'s. 

Observe that we also know $\tilde{\sigma}^2$, because it is in $G$. 
In case A) we have the expression  
$$\tilde{\sigma}^2
\left( \begin{array}{c} 
x_1 \\
y_1\\
x_2 \\
y_2
\end{array} \right)= 
\left( \begin{array}{cccc} 
1 & 0 & 0 & 0\\
0 & 1 & 0 & 0 \\
2f_1 & 0 & 1 & 0 \\
0 & 2f_1 & 0 & 1
\end{array} \right)
\left( \begin{array}{c} 
x_1 \\
y_1 \\
x_2 \\
y_2
\end{array} \right)+
\left( \begin{array}{c} 
2d_1 \\
0 \\
f_1 d_1 + 2 \gamma_1 \\
f_2 d_1 
\end{array} \right),  
$$
so we are able to determine $\gamma_1$.

Analogously in case $B)$ we find 

$$
\tilde{\sigma}^2
\left( \begin{array}{c} 
x_1 \\
y_1\\
x_2 \\
y_2
\end{array} \right)= 
\left( \begin{array}{cccc} 
1 & 0 & 0 & 0\\
0 & 1 & 0 & 0 \\
f_1 + f_2 & f_1 + f_2 & 1 & 0 \\
-(f_1 +f_2) & f_1 + f_2 & 0 & 1
\end{array} \right)
\left( \begin{array}{c} 
x_1 \\
y_1 \\
x_2 \\
y_2
\end{array} \right)+
\left( \begin{array}{c} 
0 \\
0 \\
2 \gamma_1 \\
0 
\end{array} \right).  
$$
Thus also in case B) we are able to find $\gamma_1$. 
 
Consider now a translation 
$$\tau
\left( \begin{array}{c} 
x_1 \\
y_1\\
x_2 \\
y_2
\end{array} \right) = 
\left( \begin{array}{c} 
x_1 \\
y_1 \\
x_2 + s_2\\
y_2 + t_2
\end{array} \right),$$
then 
$\tau^{-1} g \tau = g $, $\forall g \in G$, and we have

$$\tau^{-1} \tilde{\sigma} \tau
\left( \begin{array}{c} 
x_1 \\
y_1 \\
x_2 \\
y_2
\end{array} \right)
=
\left( \begin{array}{cccc} 
c_1 & c_2 &  0 & 0\\
c_2 & -c_1 & 0 & 0\\
f_1 & f_2 & 1 & 0 \\
f_2 & -f_1 & 0 & -1
\end{array} \right)
\left( \begin{array}{c} 
x_1 \\
y_1 \\
x_2 \\
y_2
\end{array} \right)+ 
\left( \begin{array}{c} 
d_1\\
0 \\
\gamma_1 \\ 
\gamma_2 - 2 t_2
\end{array} \right).
$$
If we set $t_2 = \gamma_2/2$, by substituting $\tilde{\sigma}$ with $\tau^{-1} \tilde{\sigma} \tau$ 
 we may assume that $\gamma_2 =0$ and we are done. 
\hfill{Q.E.D.} 

\begin{REM}
For a real Kodaira surface $(S, \sigma)$, we must have $\beta_2 = \frac{2iIm(\bar{\alpha_3} \alpha_4)}{m}$, $\beta_1 = Re(\beta_1) -\frac{is Im( \bar{\alpha_3} \alpha_4)}{m}$, where $Re(\beta_1) \neq 0$. Furthermore if $\tilde{\sigma}$ is any lifting of $\sigma$ to the universal covering, we have $\tilde{\sigma} g_1 \tilde{\sigma}^{-1} = g_1 g_2^s$, $\tilde{\sigma} g_2 \tilde{\sigma}^{-1} = g_2^{-1}$.

\end{REM}
{\bf Proof.}
We know that $m \beta_2 = \bar{\alpha_3} \alpha_4 - \bar{\alpha_4} \alpha_3 = 2i Im(\bar{\alpha_3} \alpha_4)$, therefore we must have $Re(\beta_1) \neq 0$, since $\beta_1$ and $\beta_2$ are linearly indipendent over ${\R}$.

Furthermore $\tilde{\sigma} g_1 \tilde{\sigma}^{-1} \in Z$, thus we have 
$$\tilde{\sigma} g_1 \tilde{\sigma}^{-1}
\left( \begin{array}{c} 
x_1\\
y_1 \\
x_2 \\
y_2
\end{array} \right)
= 
\left( \begin{array}{c} 
x_1\\
y_1 \\
x_2 \\
y_2
\end{array} \right)
+ 
\left( \begin{array}{c} 
0\\
0 \\
Re(\beta_1) \\
-Im(\beta_1)
\end{array} \right)
= g_1^r g_2^s 
\left( \begin{array}{c} 
x_1\\
y_1 \\
x_2 \\
y_2
\end{array} \right)
=$$
$$ 
\left( \begin{array}{c} 
x_1\\
y_1\\
x_2 \\
y_2
\end{array} \right)
+
\left( \begin{array}{c} 
0\\
0 \\
r Re(\beta_1) \\
rIm(\beta_1) + s Im(\beta_2)
\end{array} \right).$$
This implies $r =1$, $Im(\beta_1) = - \frac{s Im (\bar{\alpha_3} \alpha_4)}{m}$.

$$\tilde{\sigma} g_2 \tilde{\sigma}^{-1}
\left( \begin{array}{c} 
x_1\\
y_1 \\
x_2 \\
y_2
\end{array} \right)
= 
\left( \begin{array}{c} 
x_1\\
y_1 \\
x_2 \\
y_2
\end{array} \right)
+ 
\left( \begin{array}{c} 
0\\
0 \\
0 \\
-Im(\beta_2)
\end{array} \right)
= g_1^{\lambda} g_2^{\mu} 
\left( \begin{array}{c} 
x_1\\
y_1 \\
x_2 \\
y_2
\end{array} \right)
= 
$$
$$
\left( \begin{array}{c} 
x_1\\
y_1\\
x_2\\
y_2
\end{array} \right)
+
\left( \begin{array}{c} 
0\\
0 \\
\lambda Re(\beta_1) \\
\lambda Im(\beta_1) + \mu Im(\beta_2)
\end{array} \right),$$
so $\lambda =0$, $\mu =-1$.
\hfill{Q.E.D.}

\begin{REM}
\label{empty}
If $Fix (\sigma) = S({\R})\neq \emptyset$, then (\ref{OFG}) splits.
\end{REM}
{\bf Proof.}
Let $x_0$ be a fixed point of the antiholomorphic involution $\sigma$ on $S$, let $p: \tilde{S} \cong {\R}^4  \rightarrow S$ be the universal covering map. Since the covering $p$ is Galois, for all $y \in p^{-1}(x_0)$, there exists a lifting $\tilde{\sigma}$ of $\sigma$ such that $\tilde{\sigma}(y) =y$. But then $\tilde{\sigma}^2$ is a lifting of the identity map with a fixed point, therefore it must be the identity and (\ref{OFG}) splits.  
\hfill{Q.E.D.}

\begin{TEO}
Fix a topological type of a real Kodaira surface $(S, \sigma)$, then
the moduli space of the real Kodaira surfaces of the given topological
type is irreducible and connected. 
\end{TEO}
{\bf Proof.}
In section 1 we have seen that in order to give a real Kodaira surfaces it suffices to find a complex structure $\eta \in D$ on the topological type $S_0$ of $S$ such that every
$\tilde{\sigma} \in \hat{G} - G$ is antiholomorphic  the complex structure induced by $\eta$, i.e. we have $\tilde{\sigma}^* (\eta) =  \lambda \bar{\eta}$, where $\lambda \in {\Co}^*$.

We know that homogeneous coordinates on $\proj H^2(S_0, {\Co})$ are $(p_{13}, 
p_{23}, p_{14}, p_{24})$ and $\eta = p_{13} \theta_{13} + p_{23} \theta_{23} +  p_{14} \theta_{14} +  p_{24} \theta_{24} \in D$, therefore we have $p_{13} p_{24} - p_{23} p_{14} =0$, $ - p_{13} \bar{p}_{24} + p_{23} \bar{p}_{14} + p_{14} \bar{p}_{23} - p_{24} \bar{p}_{13} >0$. 

Let us choose coordinates on the universal covering of $S_0$ as above $(x_1, y_1, x_2, y_2)$. Then we have to distinguish cases A) and B) as in the proof of \ref{orfun} in the choice of a lift $\tilde{\sigma}$ of $\sigma$, $\tilde{\sigma} \in \hat{G} - G$.   

In case A) $\tilde{\sigma}$ is of the form 
$$\tilde{\sigma} 
\left( \begin{array}{c} 
x_1\\
y_1\\
x_2\\
y_2
\end{array} \right)
= 
\left( \begin{array}{cccc} 
1 & 0 & 0 & 0 \\
0 & -1 & 0 & 0\\
f_1 & f_2 & 1 & 0\\
f_2 & -f_1 & 0 & -1
\end{array} \right)
\left( \begin{array}{c} 
x_1\\
y_1\\
x_2\\
y_2
\end{array} \right)
+ 
\left( \begin{array}{c} 
d_1\\
0\\
\gamma_1\\
\gamma_2
\end{array} \right)
=
$$
$$
\left( \begin{array}{c} 
x_1 + d_1\\
-y_1\\
f_1 x_1 + f_2 y_1 + x_2 + \gamma_1\\
f_2 x_1 - f_1 y_1 - y_2 + \gamma_2 
\end{array} \right)
$$
The fact that $\tilde{\sigma}^2 \in G$ implies $f_1 = d_1$. 

Now, 
$$\tilde{\sigma}^* (\eta) = p_{13} ((f_2 - y_1)dx_1 \wedge dy_1 + dx_1 \wedge dx_2) + $$
$$p_{23} ( dx_1 \wedge dy_1 (f_1 - x_1 - d_1) - dy_1 \wedge dx_2) + $$
$$p_{14} (dx_1 \wedge dy_1(x_1 + d_1 -f_1) - dx_1 \wedge dy_2) + $$
$$p_{24} ((f_2 - y_1) dx_1 \wedge dy_1 + dy_1 \wedge dy_2) =$$
$$(dx_1 \wedge dx_2)(p_{13}) + (dx_1 \wedge dy_1)(p_{13} (f_2 - y_1) -x_1 p_{23} + x_1 p_{14} + (f_2 - y_1) p_{24}) + $$
$$(dy_1 \wedge dx_2) (-p_{23}) + (dx_1 \wedge dy_2)(-p_{14}) + (dy_1 \wedge dy_2)(p_{24}).$$ 

$$\bar{\eta} = \bar{p}_{13} (dx_1 \wedge dx_2) + (dx_1 \wedge dy_1) (-y_1 \bar{p}_{13} + x_1 \bar{p}_{23} - x_1 \bar{p}_{14} - y_1 \bar{p}_{24}) + $$
$$(dy_1 \wedge dx_2) \bar{p}_{23} + (dx_1 \wedge dy_2) \bar{p}_{14} + (dy_1 \wedge dy_2) \bar{p}_{24}.$$  

The condition $\tilde{\sigma}^*(\eta) = \lambda \bar{\eta}$ yields 
$$p_{13} = \lambda \bar{p}_{13}, \ -p_{23} = \lambda \bar{p}_{23}, \ -p_{14} = \lambda \bar{p}_{14}, \ p_{24} = \lambda \bar{p}_{24},$$
$$ p_{13}(f_2 -y_1) + p_{23}(-x_1) + p_{14} x_1 + p_{24} (f_2 - y_1) = \lambda  (-y_1 \bar{p}_{13} + x_1 \bar{p}_{23} - x_1 \bar{p}_{14} - y_1 \bar{p}_{24}),$$
so we get 
$$(p_{13} + p_{24}) f_2 = 0.$$
The point $(p_{14}/p_{24}, p_{13}/p_{14}) =:(x,y) \in (H_+ \times H_+) \cup (H_- \times H_-)$ satisfies $\bar{x} = -x$, $\bar{y} = -y$, i.e. $Re(x) = Re(y) = 0$. If $f_2=0$ these are the only conditions on $(x,y)$ therefore $(x,y) \in i {\R}_+ \times i {\R}_+ \cup i {\R}_- \times i {\R}_-$ and the two components are exchanged by the automorphism sending $(x,y)$ to $(-x, -y)$, obtained by setting in $(*)$ $a = e =-d = -1$, $b = c = k = 0$, so the moduli space is irreducible and connected. 

If $f_2 \neq 0$, we also have the condition $p_{13} = - p_{24}$, so $xy = -1$, but $x = i \gamma$, $y = i \delta$ with $(\gamma, \delta) \in {\R}_+ \times {\R}_+ \cup {\R}_- \times {\R}_-$, thus $xy = -\gamma \delta = -1$. This is a hyperbola whose two branches are exchanged by the automorphism of $D$  sending $(x,y)$ to $(-x, -y)$, therefore it is irreducible and connected.  

Now let us consider case B) for $\tilde{\sigma}$. 

$$\tilde{\sigma} 
\left( \begin{array}{c} 
x_1\\
y_1\\
x_2\\
y_2
\end{array} \right)
= 
\left( \begin{array}{cccc} 
0 & 1 & 0 & 0 \\
1 & 0 & 0 & 0\\
f_1 & f_2 & 1 & 0\\
f_2 & -f_1 & 0 & -1
\end{array} \right)
\left( \begin{array}{c} 
x_1\\
y_1\\
x_2\\
y_2
\end{array} \right)
+ 
\left( \begin{array}{c} 
0\\
0\\
\gamma_1\\
\gamma_2
\end{array} \right)
=
$$
$$
\left( \begin{array}{c} 
y_1\\
x_1\\
f_1 x_1 + f_2 y_1 + x_2 + \gamma_1\\
f_2 x_1 - f_1 y_1 - y_2 + \gamma_2
\end{array} \right)
$$
$\tilde{\sigma}^2 \in G$ implies $f_1 = -f_2$.

$$\tilde{\sigma}^*(\eta) = (dx_1 \wedge dx_2) (p_{23}) + (dx_1 \wedge dy_1)(p_{13} (x_1 -f_1) + p_{23} (f_2 - y_1) + p_{14} (y_1 -f_2) + p_{24} (x_1 -f_1))+ $$
$$(dy_1 \wedge dx_2) (p_{13}) + (dy_1 \wedge dy_2)(-p_{14}) + (dx_1 \wedge dy_2)(-p_{24}).$$
$\tilde{\sigma}^* (\eta) = \lambda \bar{\eta}$ if and only if 
$$p_{23} = \lambda \bar{p}_{13}, \ -p_{14} = \lambda \bar{p}_{24}, \ |\lambda| = 1$$
$$p_{13} (x_1 -f_1) + p_{23} (f_2 - y_1) + p_{14} (y_1- f_2) + p_{24} (x_1 -f_1) = \lambda(- y_1 \bar{p}_{13} + x_1 \bar{p}_{23} - x_1 \bar{p}_{14} - y_1 \bar{p}_{24}).$$
This last equation becomes $f_2(p_{13} + p_{24} + p_{23} - p_{14}) = 0$.

$(p_{14}/p_{24}, p_{13}/p_{14}) =: (x,y)$ satisfies $x \bar{x} = 1$, $Re(y) = \frac{1}{2p_{14}p_{24}}(-p_{14}p_{23} + p_{13} p_{24})  = 0$.

If $f_2 =0$, these are the only conditions on $(x,y) \in H_+ \times H_+ \cup H_- \times H_-$, so the moduli space is irreducible and connected, since the two components are exchanged by the same automorphism as above sending $x$ to $-x$, $y$ to $-y$. 

If $f_2 \neq 0$ we also have the condition $p_{13} + p_{23} = p_{14} - p_{24}$, that is $p_{13} + \lambda \bar{p}_{13} = p_{14} + \lambda \bar{p}_{14}$. By dividing this equation by $p_{14}$ we obtain 
$$\frac{p_{13}}{p_{14}} + \frac{\bar{p_{13}}}{\bar{p_{14}}} \frac{\lambda \bar{p_{14}}}{p_{14}} = 1 - \frac{p_{24}}{p_{14}},$$
equivalently  
$$\frac{p_{13}}{p_{14}} - \frac{\bar{p_{13}}}{\bar{p_{14}}} \frac{p_{24}}{p_{14}} = 1 - \frac{p_{24}}{p_{14}},$$
which says
$$y - \bar{y}/x = 1 - 1/x.$$
We obtain 
$$x = \frac{\bar{y} -1}{y-1} = \frac{1 +y}{1-y},$$ 
since $y + \bar{y} =0$. 
Therefore the moduli space is parametrized by the image of the map $y \mapsto (\frac{1+y}{1-y}, y)$, where $y \in H_+ \cup H_- $  and  $Re(y) = 0$. The automorphism $F$ obtained by setting in $(*)$ $a = d = 0$, $b = c =-1$, $e = ad-bc = -1$, $k =0$ sends $x = p_{14}/p_{24}$ to $1/x = \bar{x}$, $y = p_{13}/p_{14}$ to $- p_{23}/p_{24} = \bar{p_{13}}/\bar{p_{14}} = \bar{y} = -y$. So if $x =   \frac{1+y}{1-y}$, $F(x) = \bar{x} = \frac{1+\bar{y}}{1-\bar{y}}$ and the moduli space is irreducible and connected.    

\hfill{Q.E.D.}

\begin{TEO}
\label{classification}
Let us fix the topological type for a Kodaira surface $S$, i.e. we fix the torsion coefficient $m \in {\Z}$ of the first homology group of $S$. If $m \equiv 0 \ (mod \ 2)$, then the number of topologically different real Kodaira surfaces is equal to 17; if $m \equiv 1 \ (mod \ 2)$, then the number of topologically different real Kodaira surfaces is equal to 13.
\end{TEO}
{\bf Proof.}
Let us first of all fix some notation. 
We know by theorem \ref{orfun} that the topological type of a real Kodaira surface is determined by the orbifold fundamental group exact sequence
$$1 \rightarrow G \rightarrow \hat{G} \rightarrow {\Z}/2 \rightarrow 1$$
 
Every $\tilde{\sigma} \in \hat{G} - G$ has the following form:
$\tilde{\sigma} (x) = Ax +b$, where 

$$A
= \left( \begin{array}{cccc} 
c_1 & c_2 & 0 & 0\\
c_2 & -c_1 & 0 & 0\\
f_1 & f_2 & 1 & 0\\
f_2 & -f_1 & 0 & -1 
\end{array} \right),$$
$$b = 
\left( \begin{array}{c} 
d_1\\
d_2\\
\gamma_1\\
\gamma_2
\end{array} \right),
$$
while every element $h \in G$ has the form 

$$h  
\left( \begin{array}{c} 
x_1\\
y_1\\
x_2\\
y_2 
\end{array} \right)
= \left( \begin{array}{cccc} 
1 & 0 & 0 & 0\\
0 & 1 & 0 & 0\\
a & b & 1 & 0\\
-b & a & 0 & 1 
\end{array} \right)
\left( \begin{array}{c} 
x_1\\
y_1\\
x_2\\
y_2 
\end{array} \right)+ 
\left( \begin{array}{c} 
a\\
b\\
\delta\\
\epsilon
\end{array} \right).
$$
Therefore we set $g_j (x) = A_j x + v_j$, where 
$$A_j
= \left( \begin{array}{cccc} 
1 & 0 & 0 & 0\\
0 & 1 & 0 & 0\\
a_j & b_j & 1 & 0\\
-b_j & a_j & 0 & 1 
\end{array} \right),$$
$$v_j =
\left( \begin{array}{c} 
a_j\\
b_j\\
\delta_j\\
\epsilon_j
\end{array} \right),
$$
with $a_1 = a_2 =b_1 = b_2 =0$.

Since we fix the topological type of the real elliptic curve $E_{\alpha} = {\Co}/(\alpha_3 {\Z}+ \alpha_4 {\Z})$, we have for $E_{\alpha}, \sigma_1)$ three different possible topological types. Now, from the description of the moduli space of real elliptic curves (cf. e.g. \cite{ag}), it is easy to see that we may assume $a_3 =1$, $b_3 =0$, $a_4 =0$, $b_4 =1$ and we have three different topological types (cf. \ref{orfun}):

A1)  $\left( \begin{array}{cc} 
c_1 & c_2 \\
c_2 & -c_1 
\end{array} \right)= \left( \begin{array}{cc} 
1 & 0 \\
0 & -1 
\end{array} \right)$, $\left( \begin{array}{c} 
d_1 \\
d_2 
\end{array} \right)= \left( \begin{array}{c} 
0 \\
0
\end{array} \right),$ $f_1 = d_1 =0$.

A2) 
$\left( \begin{array}{cc} 
c_1 & c_2 \\
c_2 & -c_1 
\end{array} \right)= \left( \begin{array}{cc} 
1 & 0 \\
0 & -1 
\end{array} \right)$, $\left( \begin{array}{c} 
d_1 \\
d_2 
\end{array} \right)= \left( \begin{array}{c} 
1/2 \\
0
\end{array} \right),$ $f_1 = d_1 = 1/2$.

B) $\left( \begin{array}{cc} 
c_1 & c_2 \\
c_2 & -c_1 
\end{array} \right)= \left( \begin{array}{cc} 
0 & 1 \\
1 & 0 
\end{array} \right)$, $\left( \begin{array}{c} 
d_1 \\
d_2 
\end{array} \right)= \left( \begin{array}{c} 
0 \\
0
\end{array} \right),$ $f_2 = -f_1$. 

Finally, as we showed in the proof of theorem \ref{orfun}, by conjugating by some translation $\tau \left( \begin{array}{c} 
x_1 \\
y_1 \\
x_2 \\
y_2 
\end{array} \right)= \left( \begin{array}{c} 
x_1 \\
y_1 \\
x_2 + s_2 \\
y_2 + t_2
\end{array} \right)$, we have $\tau^{-1} g \tau = g$, for all $g \in G$ and $\tau^{-1} \tilde{\sigma} \tau (x) = Ax + \left( \begin{array}{c} 
d_1 \\
0 \\
\gamma_1\\
\gamma_2 - 2t_2
\end{array} \right)$ and so we may assume that in the expression of $\tilde{\sigma}$ we have $\gamma_2 =0$.   

Recall that $g_3 g_4 g_3^{-1} g_4^{-1} = g_2^m$, so $\delta_2 =0$ and $\epsilon_2 = 2/m$, and thus $\delta_1 \neq 0$. 

Observe that for every $g \in Z$, $g(x) = x +v$, and for every $\tilde{\sigma} \in \hat{G} - G$ we have 
$\tilde{\sigma} g \tilde{\sigma}^{-1}(x) = x + Av$, therefore we have 
  
$\tilde{\sigma} g_2 \tilde{\sigma}^{-1}(x) = x + Av_2 = g_2^{-1}(x)$, since $\delta_2 =0$. 

$\tilde{\sigma} g_1 \tilde{\sigma}^{-1}(x) = x + Av_1 = x + \left( \begin{array}{c} 
0 \\
0 \\
\delta_1 \\
-\epsilon_1
\end{array} \right)$. Since we know that $\tilde{\sigma} Z \tilde{\sigma}^{-1} = Z$, we must have $\tilde{\sigma} g_1 \tilde{\sigma}^{-1} = g_1^{\lambda} g_2^{\mu}$ and we immediately see that we have $\lambda = 1$ $\mu = -2 \epsilon_1/\epsilon_2 = -m \epsilon_1 \in {\Z}$.

So we have shown that for every $\tilde{\sigma} \in \hat{G} - G$ we have $\tilde{\sigma} g_2 \tilde{\sigma}^{-1} = g_2^{-1}$, $\tilde{\sigma} g_1 \tilde{\sigma}^{-1} = g_1 g_2^{\mu}$, $\mu = -\frac{2\epsilon_1}{\epsilon_2} = -m \epsilon_1 \in {\Z}$. 

Recall that we have an action of ${\Z}/2 =<\sigma>$ on $Z$ given by the orbifold fundamental group exact sequence, since the conjugation action on $Z$ is independent of the choice of the lifting $\tilde{\sigma} \in \hat{G} - G$ of $\sigma$. \\
 
Observe that we still can change generators of $Z$ by substituting $g_1$ with $g_1 g_2^t$, and leaving $g_2$ invariant. 

Then we have $\tilde{\sigma} g_1 g_2^t \tilde{\sigma}^{-1} = g_1 g_2^{\mu} g_2^{-t}$. Thus if $\mu \equiv 0 \ (mod \ 2)$, by choosing $t = \mu/2$ we can assume $\tilde{\sigma} g_1 \tilde{\sigma}^{-1} = g_1$, $\mu = \epsilon_1 =0$, while if $\mu \equiv 1 \ (mod \ 2)$, by choosing $t = (\mu-1)/2$ we can assume $\tilde{\sigma} g_1 \tilde{\sigma}^{-1} = g_1 g_2$, $\mu =1$, $\epsilon_1 = -1/m$. 

So for every $\tilde{\sigma} \in \hat{G} - G$, for the action of $\tilde{\sigma}$ on $Z$ we have the two following different cases:\\

1) $\tilde{\sigma} g_1 \tilde{\sigma}^{-1} = g_1$, $\mu = \epsilon_1 =0$, $\tilde{\sigma} g_2 \tilde{\sigma}^{-1} = g_2^{-1},$\\

2) $\tilde{\sigma} g_1 \tilde{\sigma}^{-1} = g_1 g_2$, $\mu = 1$,  $\epsilon_1 =-1/m$, $\tilde{\sigma} g_2 \tilde{\sigma}^{-1} = g_2^{-1}.$\\

Observe that since we have fixed the action of a lifting $\tilde{\sigma}$ of $\sigma$ on the elliptic curve $E_{\alpha} = {\Co}/(\alpha_3 {\Z} + \alpha_4{\Z})$, we can now change lifting only by multiplying $\tilde{\sigma}$ by an element $g \in Z$. 
In fact if we take a lifting $\tilde{\sigma}' = \tilde{\sigma} g$ with $g = g_4^b g_3^a g_1^r g_2^s$ we see that the linear part of the action of $\tilde{\sigma}'$ on the first two coordinates does not change, while the translation part changes. 
In fact we always can assume, by conjugating by a translation, that in the new lifting we have $\gamma'_2 =0$. Furthermore we see that if $\tilde{\sigma}$ is chosen as in cases A1) and A2) we obtain for $\tilde{\sigma}' = \tilde{\sigma} g$, $d'_1 = d_1 +a$, $d'_2 = d_2 -b = -b$, while in case B) we have $d'_1 = d_1 +b$, $d'_2 = d_2 + a$.

Let us then consider the action of an element $\tilde{\sigma} \in \hat{G} - G$ on $G - Z$. \\

In case B) we easily see that we have 
$$\tilde{\sigma} g_3 \tilde{\sigma}^{-1}
\left( \begin{array}{c} 
x_1\\
y_1\\
x_2\\
y_2 
\end{array} \right)
= \left( \begin{array}{cccc} 
1 & 0 & 0 & 0\\
0 & 1 & 0 & 0\\
0 & 1 & 1 & 0\\
-1 & 0 & 0 & 1 
\end{array} \right)
\left( \begin{array}{c} 
x_1\\
y_1\\
x_2\\
y_2 
\end{array} \right)+ 
\left( \begin{array}{c} 
0\\
1\\
f_1 + \delta_3\\
-f_1 - \epsilon_3
\end{array} \right),
$$
therefore, looking at the linear parts we conclude that we must have $\tilde{\sigma} g_3 \tilde{\sigma}^{-1} = g_4 g_1^r g_2^n$.\\

Observe that by substituting $g_3$ with $g_3 g_1^l g_2^t$ we do not change the topological type of the Kodaira surface and we obtain $\tilde{\sigma} (g_3 g_1^l g_2^t) \tilde{\sigma}^{-1} = g_4 g_1^{r+l} g_2^{n + l \mu -t}$. 

Thus if we choose $l = -r$, $t = n - \mu r$, we see that we can assume $r =0$, $n = 0$, $\tilde{\sigma} g_3 \tilde{\sigma}^{-1} = g_4$. Then we obtain  $f_1 = \delta_4 -\delta_3 = -(\epsilon_4 + \epsilon_3)$. \\
By an easy computation one shows that $\tilde{\sigma} g_4 \tilde{\sigma}^{-1} = g_3$.\\

In cases A1) and A2) we have 
$$\tilde{\sigma} g_3 \tilde{\sigma}^{-1}
\left( \begin{array}{c} 
x_1\\
y_1\\
x_2\\
y_2 
\end{array} \right)
= \left( \begin{array}{cccc} 
1 & 0 & 0 & 0\\
0 & 1 & 0 & 0\\
1 & 0 & 1 & 0\\
0 & 1 & 0 & 1 
\end{array} \right)
\left( \begin{array}{c} 
x_1\\
y_1\\
x_2\\
y_2 
\end{array} \right)+ 
\left( \begin{array}{c} 
1\\
0\\
\delta_3\\
f_2 -\epsilon_3
\end{array} \right),
$$
so by looking at the linear parts we see that we must have 
$$\tilde{\sigma} g_3 \tilde{\sigma}^{-1}
\left( \begin{array}{c} 
x_1\\
y_1\\
x_2\\
y_2 
\end{array} \right)
= g_3 g_1^r g_2^n \left( \begin{array}{c} 
x_1\\
y_1\\
x_2\\
y_2 
\end{array} \right)
= \left( \begin{array}{cccc} 
1 & 0 & 0 & 0\\
0 & 1 & 0 & 0\\
1 & 0 & 1 & 0\\
0 & 1 & 0 & 1 
\end{array} \right)
\left( \begin{array}{c} 
x_1\\
y_1\\
x_2\\
y_2 
\end{array} \right)+ 
\left( \begin{array}{c} 
1\\
0\\
r \delta_1 + \delta_3\\
r \epsilon_1 + n \epsilon_2 + \epsilon_3
\end{array} \right).
$$
Thus we obtain $r \delta_1 = 0$, which implies $ r =0$, since $\delta_1 \neq 0$, $f_2 = 2 \epsilon_3 + n \epsilon_2$.
So we have $\tilde{\sigma} g_3 \tilde{\sigma}^{-1} = g_3 g_2^n$.  

By substituting $g_3$ with $g_3 g_1^l g_2^t$ we do not change the topological type of the Kodaira surface and we have
$\tilde{\sigma} g_3 g_1^l g_2^t \tilde{\sigma}^{-1} = g_3 g_2^n g_1^l g_2^{\mu l} g_2^{-t} = g_3 g_1^l g_2^{n + \mu l -t}$.\\

In case 1) we have $\mu = \epsilon_1 =0$ and therefore if $n \equiv 0$ $(mod \ 2)$ we may choose $t = n/2$, $l=0$ and by substituting $g_3$ with $g_3 g_2^{n/2}$, we see that we can assume $n =0$, $\tilde{\sigma} g_3 \tilde{\sigma}^{-1} = g_3$, and $f_2 = 2 \epsilon_3$. 

If $n   \equiv 1$ $(mod \ 2)$, we can set $t = (n -1)/2$, $l=0$ and by substituting $g_3$ with $g_3 g_2^{(n-1)/2}$, we see that we can assume $n = 1$, $\tilde{\sigma} g_3 \tilde{\sigma}^{-1} = g_3 g_2$, and $f_2 = 2 \epsilon_3 + \epsilon_2$.\\

In case 2) we have $\mu = 1$, $ \epsilon_1 =-1/m$ and therefore if $n \equiv 0$ $(mod \ 2)$ we may choose $l =0$, $t = n/2$ and by substituting $g_3$ with $g_3 g_2^{n/2}$, we see that we can assume $n=0$, $\tilde{\sigma} g_3 \tilde{\sigma}^{-1} = g_3$, and $f_2 = 2 \epsilon_3$.

If $n \equiv 1$ $(mod \ 2)$, we may choose $l =1$, $t = (n +1)/2$ and by substituting $g_3$ with $g_3 g_1 g_2^{(n + 1)/2}$, we see that we can assume $n =0$, $\tilde{\sigma} g_3 \tilde{\sigma}^{-1} = g_3$, and $f_2 = 2 \epsilon_3$. \\

In conclusion we have the following cases:\\

1)A1)A2)(a): $\tilde{\sigma} g_1 \tilde{\sigma}^{-1} = g_1$, $\tilde{\sigma} g_2 \tilde{\sigma}^{-1} = g_2^{-1}$, ($\mu = \epsilon_1 =0$), $\tilde{\sigma} g_3 \tilde{\sigma}^{-1} = g_3$, ($f_2 = 2 \epsilon_3$). \\

1)A1)A2)(b): $\tilde{\sigma} g_1 \tilde{\sigma}^{-1} = g_1$, $\tilde{\sigma} g_2 \tilde{\sigma}^{-1} = g_2^{-1}$, ($\mu = \epsilon_1 =0$), $\tilde{\sigma} g_3 \tilde{\sigma}^{-1} = g_2 g_3$, ($f_2 = 2 \epsilon_3 + \epsilon_2$).\\

2)A1)A2): $\tilde{\sigma} g_1 \tilde{\sigma}^{-1} = g_1 g_2$, $\tilde{\sigma} g_2 \tilde{\sigma}^{-1} = g_2^{-1}$, ($\mu = 1$, $\epsilon_1 =-1/m$), $\tilde{\sigma} g_3 \tilde{\sigma}^{-1} = g_3$,  ($f_2 = 2 \epsilon_3$). \\

By an easy computation we see that in cases A1) and A2) we have 

$$\tilde{\sigma} g_4 \tilde{\sigma}^{-1} = \left( \begin{array}{c} 
x_1\\
y_1\\
x_2\\
y_2 
\end{array} \right)
\left( \begin{array}{cccc} 
1 & 0 & 0 & 0\\
0 & 1 & 0 & 0\\
0 & -1 & 1 & 0\\
1 & 0 & 0 & 1 
\end{array} \right)
\left( \begin{array}{c} 
x_1\\
y_1\\
x_2\\
y_2 
\end{array} \right)+ 
\left( \begin{array}{c} 
0\\
-1\\
f_2 +\delta_4\\
-2d_1 - \epsilon_4
\end{array} \right),
$$

so by looking at the linear parts we conclude that we must have $\tilde{\sigma} g_4 \tilde{\sigma}^{-1} = g_4^{-1} g_1^u g_2^v$. Now we have

$$g_4^{-1} g_1^u g_2^v \left( \begin{array}{c} 
x_1\\
y_1\\
x_2\\
y_2 
\end{array} \right)
= \left( \begin{array}{cccc} 
1 & 0 & 0 & 0\\
0 & 1 & 0 & 0\\
0 & -1 & 1 & 0\\
1 & 0 & 0 & 1 
\end{array} \right)
\left( \begin{array}{c} 
x_1\\
y_1\\
x_2\\
y_2 
\end{array} \right)+ 
\left( \begin{array}{c} 
0\\
-1\\
1 - \delta_4 +u \delta_1\\
- \epsilon_4 + u \epsilon_1 + v \epsilon_2
\end{array} \right),$$
thus we obtain $u = (f_2 + 2 \delta_4 -1)/\delta_1$, $v =  (-2d_1 -u \epsilon_1)/\epsilon_2$.

Again we can substitute $g_4$ with $g_4 g_1^r g_2^s$ and try to find a normal form for the conjugation action of $\tilde{\sigma}$ on $g_4$. 

In case 1)A1)(a) we compute $u = (2 \epsilon_3 -1 + 2 \delta_4)/\delta_1 \in {\Z}$, $v = (-2d_1 -u \epsilon_1)/\epsilon_2 =0$, $\tilde{\sigma} g_4 g_1^r g_2^s \tilde{\sigma}^{-1} = g_4^{-1} g_1^u  g_1^r g_2^{-s}$.\\
Thus we have the two following possibilities:\\

1)A1)(a)(i): If $u \equiv 0$ $(mod \ 2)$, then we can set $r = -u/2$ and we may assume $\tilde{\sigma} g_4 \tilde{\sigma}^{-1} = g_4^{-1}$, $u =0$, i.e. $\epsilon_3 + \delta_4 =1/2$.\\

1)A1)(a)(ii): If $u \equiv 1 \ (mod \ 2)$, then we can set $r = (1-u)/2$ and we may assume $\tilde{\sigma} g_4 \tilde{\sigma}^{-1} = g_4^{-1}g_1$, $u =1$, i.e. $2\epsilon_3 + 2\delta_4 -1=\delta_1$.\\

In case 1)A1)(b) we compute $u = (2 \epsilon_3 -1 + 2 \delta_4 +\epsilon_2)/\delta_1 \in {\Z}$, $v = (-2d_1 -u \epsilon_1)/\epsilon_2 =0$, $\tilde{\sigma} g_4 g_1^r g_2^s \tilde{\sigma}^{-1} = g_4^{-1} g_1^u g_1^r g_2^{-s}$.\\
Thus we have the two following possibilities:\\

1)A1)(b)(i): If $u \equiv 0 \ (mod \ 2)$, then we can set $r = -u/2$ and we may assume $\tilde{\sigma} g_4 \tilde{\sigma}^{-1} = g_4^{-1}$, $u =0$, i.e. $2\epsilon_3 + 2\delta_4+ \epsilon_2 =1$.\\

1)A1)(b)(ii): If $u \equiv 1 \ (mod \ 2)$, then we can set $r = (1-u)/2$ and we may assume $\tilde{\sigma} g_4 \tilde{\sigma}^{-1} = g_4^{-1}g_1$, $u =1$, i.e. $2\epsilon_3 + 2\delta_4 -1 + \epsilon_2 =\delta_1$.\\

In case 1)A2)(a) we compute $u = (2 \epsilon_3 -1 + 2 \delta_4)/\delta_1 \in {\Z}$, $v =(-2d_1 -u \epsilon_1)/\epsilon_2 =-m/2$, therefore this case occurs only if $m \equiv 0 \ (mod \ 2)$. $\tilde{\sigma} g_4 g_1^r g_2^s \tilde{\sigma}^{-1} = g_4^{-1} g_1^u g_2^{\frac{-m}{2}}g_1^r g_2^{-s}$.\\
Thus we have the two following possibilities:\\

1)A2)(a)(i): If $u \equiv 0 \ (mod \ 2)$, then we can set $r = -u/2$ and we may assume $\tilde{\sigma} g_4 \tilde{\sigma}^{-1} = g_4^{-1} g_2^{-m/2}$, $u =0$, i.e. $\epsilon_3 + \delta_4 =1/2$.\\

1)A2)(a)(ii): If $u \equiv 1 \ (mod \ 2)$, then we can set $r = (1-u)/2$ and we may assume $\tilde{\sigma} g_4 \tilde{\sigma}^{-1} = g_4^{-1}g_1 g_2^{-m/2}$, $u =1$, i.e. $2\epsilon_3 + 2\delta_4 -1=\delta_1$.\\

We will prove later that case 1)A2)(b) does not occur.\\

In case 2)A1) we have $\tilde{\sigma} g_2 \tilde{\sigma}^{-1} = g_2^{-1}$, $\tilde{\sigma} g_1 \tilde{\sigma}^{-1} = g_1 g_2$, ($\mu = 1$, $\epsilon_1 = -1/m$), $\tilde{\sigma} g_3 \tilde{\sigma}^{-1} = g_3$ ($f_2 = 2 \epsilon_3$), $\tilde{\sigma} g_4 \tilde{\sigma}^{-1} = g_4^{-1} g_1^u g_2^v$, with $u = (2 \epsilon_3 + 2\delta_4 -1)/\delta_1$, $v =(-2d_1 -u \epsilon_1)/\epsilon_2 =  u/2$. Thus we must have $u \equiv 0 \ (mod \ 2)$. \\
Then $\tilde{\sigma} g_4 g_1^r g_2^s \tilde{\sigma}^{-1} = g_4^{-1} g_1^u g_2^v g_1^r g_2^{r-s} $, so if we choose $r = -u/2$ we have $v + r = 0$ and we can assume $\tilde{\sigma} g_4 \tilde{\sigma}^{-1} = g_4^{-1}$, $u =0$, $v = 0$, $\epsilon_3 + \delta_4 = 1/2$. \\

In case 2)A2) we have $\tilde{\sigma} g_2 \tilde{\sigma}^{-1} = g_2^{-1}$, $\tilde{\sigma} g_1 \tilde{\sigma}^{-1} = g_1 g_2$, ($\mu = 1$, $\epsilon_1 = -1/m$), $\tilde{\sigma} g_3 \tilde{\sigma}^{-1} = g_3$ ($f_2 = 2 \epsilon_3$), $\tilde{\sigma} g_4 \tilde{\sigma}^{-1} = g_4^{-1} g_1^u g_2^v$, with $u = (2 \epsilon_3 + 2\delta_4 -1 )/\delta_1$, $v = (-2d_1 -u \epsilon_1)/\epsilon_2 = (u-m)/2$. Thus we must have $u \equiv m \ (mod \ 2)$. \\
Then $\tilde{\sigma} g_4 g_1^r g_2^s \tilde{\sigma}^{-1} = g_4^{-1} g_1^u g_2^v g_1^r g_2^{r-s} $, so we have the following two different cases:\\

2)A2)(i) If $u \equiv 0 \ (mod \ 2)$ (then also $m \equiv 0 \ (mod \ 2)$), we can choose $r = -u/2$ and we have $v + r =  -m/2$ and we can assume $\tilde{\sigma} g_4 \tilde{\sigma}^{-1} = g_4^{-1}g_2^{-m/2}$, $v = -m/2$, $u = 0$, i.e. $\epsilon_3 + \delta_4 = 1/2$. \\

2)A2)(ii) If $u \equiv 1 \ (mod \ 2)$ (then also $m \equiv 1 \ (mod \ 2)$), we can choose $r = (-u+1)/2$ and we have $v + r =  (1-m)/2$ and we can assume $\tilde{\sigma} g_4 \tilde{\sigma}^{-1} = g_4^{-1}g_1 g_2^{(1-m)/2}$, $v =  (1-m)/2$, $u = 1$, i.e. $2\epsilon_3 + 2\delta_4 -1 = \delta_1$. \\

Observe now that in the cases A1) and B) we have $\tilde{\sigma}^2 \in Z$, $\forall \tilde{\sigma} \in \hat{G} - G$. In fact,  

$$\tilde{\sigma}^2
\left( \begin{array}{c} 
x_1\\
y_1\\
x_2\\
y_2 
\end{array} \right)
= \left( \begin{array}{cccc} 
1 & 0 & 0 & 0\\
0 & 1 & 0 & 0\\
0 & 0 & 1 & 0\\
0 & 0 & 0 & 1 
\end{array} \right)
\left( \begin{array}{c} 
x_1\\
y_1\\
x_2\\
y_2 
\end{array} \right)+ 
\left( \begin{array}{c} 
0\\
0\\
2\gamma_1 \\
0
\end{array} \right).
$$
So, in cases A1) and B) we can write $\tilde{\sigma}^2 = g_1^p g_2^q \in Z$. 

In case A2), for every lifting $\tilde{\sigma}$, we have 
$$\tilde{\sigma}^2
\left( \begin{array}{c} 
x_1\\
y_1\\
x_2\\
y_2 
\end{array} \right)
= \left( \begin{array}{cccc} 
1 & 0 & 0 & 0\\
0 & 1 & 0 & 0\\
1 & 0 & 1 & 0\\
0 & 1 & 0 & 1 
\end{array} \right)
\left( \begin{array}{c} 
x_1\\
y_1\\
x_2\\
y_2 
\end{array} \right)+ 
\left( \begin{array}{c} 
1\\
0\\
1/4 + 2\gamma_1 \\
f_2/2
\end{array} \right),
$$
so we can write $\tilde{\sigma}^2 = g_3 g_1^p g_2^q $, thus (\ref{OFG}) doesn't split.

Therefore in cases A1) and B) we must have $2 \gamma_1 = p \delta_1$, $p \epsilon_1 + q \epsilon_2 =0$, or equivalently $p \mu = 2 q$. 

In case A2) we must have $2 \gamma_1  = -1/4 + p \delta_1 + \delta_3$, $p \epsilon_1 + q \epsilon_2 + \epsilon_3 = f_2/2$. \\

Notice that by substituting $\tilde{\sigma}$ by $\tilde{\sigma} g_1^k g_2^s$, $k, s \in {\Z}$, the conjugation action of $\tilde{\sigma}$ on $G$ doesn't change. 
So if we substitute $\tilde{\sigma}$ by $\tilde{\sigma} g_1^k g_2^s$, in cases A1) and B) we get $(\tilde{\sigma} g_1^k g_2^s)^2 = g_1^{2k + p} g_2^{k \mu  +q}$.

In case 1)B) then we have $\tilde{\sigma}^2 = g_1^p g_2^q$, $2 \gamma_1 = p \delta_1$, $p \epsilon_1 + q \epsilon_2 = q \epsilon_2 =0$, thus $q =0$, since $\epsilon_2 \neq 0$. So if we substitute $\tilde{\sigma}$ by $\tilde{\sigma} g_1^k g_2^s$, we obtain 
$(\tilde{\sigma} g_1^k g_2^s)^2 = g_1^{2k +p}$, so we have the two following cases:\\

1)B)' if $p \equiv 0 \ (mod \ 2)$, then we can set $k = -p/2$ and we may assume $\tilde{\sigma}^2 = Id$, $\gamma_1 =0$ and (\ref{OFG}) splits.\\

1)B)'' if $p \equiv 1 \ (mod \ 2)$, then we can set $k = (1-p)/2$ and we may assume $\tilde{\sigma}^2 = g_1$, $\gamma_1 =\delta_1/2$. We claim that (\ref{OFG}) does not split. In fact otherwise there would exist an element $g \in G$ such that $(\tilde{\sigma} g)^2 =1$. But any $g \in G$ can be written as follows: $g = g_4^r g_3^s g_1^l g_2^t$, so we have $(\tilde{\sigma}g)^2 = (\tilde{\sigma} g \tilde{\sigma}^{-1}) \tilde{\sigma} (\tilde{\sigma} g \tilde{\sigma}^{-1}) \tilde{\sigma} = g_3^r g_4^s g_1^l g_2^{-t} \tilde{\sigma}(g_3^r g_4^s g_1^l g_2^{-t}) \tilde{\sigma}^{-1} \tilde{\sigma}^2= g_3^r g_4^{s+r} g_3^s g_1^{2l+1} = g_4^{s+r} g_3^{s+r} g_2^{mr(r+s)}g_1^{2l+1}$. Then $(\tilde{\sigma}g)^2 = 1$ implies $2l +1 =0$, a contradiction. \\

In case 2)B) we have $\tilde{\sigma}^2 = g_1^p g_2^q$, $2 \gamma_1 = p \delta_1$, $p \epsilon_1 + q \epsilon_2 = \frac{-p +2q}{m} =0$, thus $2q =p$. So if we substitute $\tilde{\sigma}$ by $\tilde{\sigma} g_1^k g_2^s$, we obtain 
$(\tilde{\sigma} g_1^k g_2^s)^2 = g_1^{2k +p}g_2^{k+q} = g_1^{2(k +q)}g_2^{k+q} $. Therefore if we set $k = -q$ we can assume $\tilde{\sigma}^2 =Id$, $\gamma_1 =0$ and (\ref{OFG}) splits.\\

In case 1)A1) we have $\tilde{\sigma}^2 = g_1^p g_2^q$, $2 \gamma_1 = p \delta_1$, $p \epsilon_1 + q \epsilon_2 = q \epsilon_2 =0$, thus $q =0$, since $\epsilon_2 \neq 0$. So if we substitute $\tilde{\sigma}$ by $\tilde{\sigma} g_1^k g_2^s$, we obtain 
$(\tilde{\sigma} g_1^k g_2^s)^2 = g_1^{2k +p}$, so we have the two following cases:\\

1)A1)' if $p \equiv 0 \ (mod \ 2)$, then we can set $k = -p/2$ and we may assume $\tilde{\sigma}^2 = Id$, $\gamma_1 =0$ and (\ref{OFG}) splits.\\

1)A1)'' if $p \equiv 1 \ (mod \ 2)$, then we can set $k = (1-p)/2$ and we may assume $\tilde{\sigma}^2 = g_1$, $\gamma_1 =\delta_1/2$. Observe that this does not necessarily imply that   (\ref{OFG}) does not split.\\

In case 2)A1) we have $\tilde{\sigma}^2 = g_1^p g_2^q$, $2 \gamma_1 = p \delta_1$, $p \epsilon_1 + q \epsilon_2 = \frac{-p +2q}{m} =0$, thus $2q =p$. So if we substitute $\tilde{\sigma}$ by $\tilde{\sigma} g_1^k g_2^s$, we obtain 
$(\tilde{\sigma} g_1^k g_2^s)^2 = g_1^{2k +p}g_2^{k+q} = g_1^{2(k +q)}g_2^{k+q} $. Therefore if we set $k = -q$ we can assume $\tilde{\sigma}^2 =Id$, $\gamma_1 =0$ and (\ref{OFG}) splits.\\

In case 1)A2)(a) we have $\tilde{\sigma}^2 = g_3g_1^p g_2^q$, $2 \gamma_1 = -1/4 + p \delta_1 +\delta_3$, $p \epsilon_1 + q \epsilon_2 + \epsilon_3 = q \epsilon_2 + \epsilon_3= f_2/2 = \epsilon_3$, thus $q =0$, since $\epsilon_2 \neq 0$. So if we substitute $\tilde{\sigma}$ by $\tilde{\sigma} g_1^k g_2^s$, we obtain 
$(\tilde{\sigma} g_1^k g_2^s)^2 = g_3 g_1^{2k +p}$, so we have the two following cases:\\

1)A2)(a)' if $p \equiv 0 \ (mod \ 2)$, then we can set $k = -p/2$ and we may assume $\tilde{\sigma}^2 = g_3$, $\gamma_1 =(\delta_3-1/4)/2$.\\

1)A2)(a)'' if $p \equiv 1 \ (mod \ 2)$, then we can set $k = (1-p)/2$ and we may assume $\tilde{\sigma}^2 = g_3 g_1$, $\gamma_1 =(-1/4 + \delta_1 + \delta_3)/2$.\\

In case 1)A2)(b) we have $2 \gamma_1 = -1/4 + p \delta_1 +\delta_3$, $p \epsilon_1 + q \epsilon_2 + \epsilon_3 = q \epsilon_2 + \epsilon_3 = f_2/2 = \epsilon_3 + \epsilon_2/2$, thus $q =1/2$, absurd. Thus this case does not occur. \\

In case 2)A2) we have $\tilde{\sigma}^2 = g_3g_1^p g_2^q$, $2 \gamma_1 = -1/4 + p \delta_1 +\delta_3$, $p \epsilon_1 + q \epsilon_2 + \epsilon_3 = -p/m + 2q/m + \epsilon_3= f_2/2 = \epsilon_3$, thus $2q =p$. So if we substitute $\tilde{\sigma}$ by $\tilde{\sigma} g_1^k g_2^s$, we obtain 
$(\tilde{\sigma} g_1^k g_2^s)^2 = g_3 g_1^{2k +p}g_2^{k +q} = g_3 g_1^{2(k +q)}g_2^{k +q}$, so if we set $k =-q$, we can assume $\tilde{\sigma}^2 = g_3$, $\gamma_1 = (\delta_3 -1/4)/2$. \\

We finally list all the different cases in normal forms that occur:\\

1)B)' $(E_{\alpha}, \sigma_1)$ as in case B), $\gamma_2 =0$, $\epsilon_1 =0$. $\tilde{\sigma} g_2 \tilde{\sigma}^{-1} = g_2^{-1}$, $\tilde{\sigma} g_1 \tilde{\sigma}^{-1} = g_1$, $\tilde{\sigma} g_3 \tilde{\sigma}^{-1} = g_4$, $f_1 = -f_2 = \delta_4 -\delta_3 =  -\epsilon_4 -\epsilon_3$, $\tilde{\sigma} g_4 \tilde{\sigma}^{-1} = g_3$, $\tilde{\sigma}^2 = Id$, $\gamma_1 =0$ and (\ref{OFG}) splits.\\

1)B)'' $(E_{\alpha}, \sigma_1)$ as in case B), $\gamma_2 =0$, $\epsilon_1 =0$. $\tilde{\sigma} g_2 \tilde{\sigma}^{-1} = g_2^{-1}$, $\tilde{\sigma} g_1 \tilde{\sigma}^{-1} = g_1$, $\tilde{\sigma} g_3 \tilde{\sigma}^{-1} = g_4$, $f_1 = -f_2 = \delta_4 -\delta_3 =  -\epsilon_4 -\epsilon_3$, $\tilde{\sigma} g_4 \tilde{\sigma}^{-1} = g_3$, $\tilde{\sigma}^2 = g_1$, $\gamma_1 =\delta_1/2$ and we have already observed that(\ref{OFG}) does not split.\\

2)B) $(E_{\alpha}, \sigma_1)$ as in case B), $\gamma_2 =0$, $\epsilon_1 =-1/m$.
$\tilde{\sigma} g_2 \tilde{\sigma}^{-1} = g_2^{-1}$, $\tilde{\sigma} g_1 \tilde{\sigma}^{-1} = g_1 g_2$, $\tilde{\sigma} g_3 \tilde{\sigma}^{-1} = g_4$, $f_1 = \delta_4 -\delta_3 = -\epsilon_4 -\epsilon_3$, $\tilde{\sigma} g_4 \tilde{\sigma}^{-1} = g_3$, $\tilde{\sigma}^2 = Id$, $\gamma_1 =0$ and (\ref{OFG}) splits.\\

1)A1)(a)(i)' $(E_{\alpha}, \sigma_1)$ as in case A1), $\gamma_2 =0$, $\epsilon_1 =0$. $\tilde{\sigma} g_2 \tilde{\sigma}^{-1} = g_2^{-1}$, $\tilde{\sigma} g_1 \tilde{\sigma}^{-1} = g_1$, $\tilde{\sigma} g_3 \tilde{\sigma}^{-1} = g_3$, $f_2 = 2 \epsilon_3$, $\tilde{\sigma} g_4 \tilde{\sigma}^{-1} = g_4^{-1}$, $\epsilon_3 + \delta_4 = 1/2$, $\tilde{\sigma}^2 = Id$, $\gamma_1 =0$ and (\ref{OFG}) splits.\\

1)A1)(a)(i)'' $(E_{\alpha}, \sigma_1)$ as in case A1), $\gamma_2 =0$, $\epsilon_1 =0$. $\tilde{\sigma} g_2 \tilde{\sigma}^{-1} = g_2^{-1}$, $\tilde{\sigma} g_1 \tilde{\sigma}^{-1} = g_1$, $\tilde{\sigma} g_3 \tilde{\sigma}^{-1} = g_3$, $f_2 = 2 \epsilon_3$, $\tilde{\sigma} g_4 \tilde{\sigma}^{-1} = g_4^{-1}$, $\epsilon_3 + \delta_4 = 1/2$, $\tilde{\sigma}^2 = g_1$, $\gamma_1 =\delta_1/2$. We also claim that (\ref{OFG}) does not split. In fact otherwise there would exist an element $g \in G$ such that $(\tilde{\sigma} g)^2 =1$. But any $g \in G$ can be written as follows: $g = g_4^r g_3^s g_1^l g_2^t$, so we have $(\tilde{\sigma}g)^2 = (\tilde{\sigma} g \tilde{\sigma}^{-1}) \tilde{\sigma} (\tilde{\sigma} g \tilde{\sigma}^{-1}) \tilde{\sigma} = g_4^{-r} g_3^s g_1^l g_2^{-t} \tilde{\sigma}(g_4^{-r} g_3^s g_1^l g_2^{-t}) \tilde{\sigma}^{-1} \tilde{\sigma}^2= g_4^{-r} g_3^{s} g_4^r g_3^s g_1^{2l+1} = g_3^{2s} g_2^{mrs} g_1^{2l+1}$. Then $(\tilde{\sigma}g)^2 = 1$ implies $2l +1 =0$, a contradiction. \\

1)A1)(a)(ii)' $(E_{\alpha}, \sigma_1)$ as in case A1), $\gamma_2 =0$, $\epsilon_1 =0$. $\tilde{\sigma} g_2 \tilde{\sigma}^{-1} = g_2^{-1}$, $\tilde{\sigma} g_1 \tilde{\sigma}^{-1} = g_1$, $\tilde{\sigma} g_3 \tilde{\sigma}^{-1} = g_3$, $f_2 = 2 \epsilon_3$, $\tilde{\sigma} g_4 \tilde{\sigma}^{-1} = g_4^{-1}g_1$, $2\epsilon_3 + 2\delta_4 -1= \delta_1$, $\tilde{\sigma}^2 = Id$, $\gamma_1 =0$ and (\ref{OFG}) splits.\\

1)A1)(a)(ii)'' $(E_{\alpha}, \sigma_1)$ as in case A1), $\gamma_2 =0$, $\epsilon_1 =0$. $\tilde{\sigma} g_2 \tilde{\sigma}^{-1} = g_2^{-1}$, $\tilde{\sigma} g_1 \tilde{\sigma}^{-1} = g_1$, $\tilde{\sigma} g_3 \tilde{\sigma}^{-1} = g_3$, $f_2 = 2 \epsilon_3$, $\tilde{\sigma} g_4 \tilde{\sigma}^{-1} = g_4^{-1}g_1$, $2\epsilon_3 + 2\delta_4 -1= \delta_1$, $\tilde{\sigma}^2 = g_1$, $\gamma_1 =\delta_1/2$ and we claim that (\ref{OFG}) splits. In fact one can easily see that $(\tilde{\sigma} g_4 g_1^{-1})^2 = 1$.\\

1)A1)(b)(i)' $(E_{\alpha}, \sigma_1)$ as in case A1), $\gamma_2 =0$, $\epsilon_1 =0$. $\tilde{\sigma} g_2 \tilde{\sigma}^{-1} = g_2^{-1}$, $\tilde{\sigma} g_1 \tilde{\sigma}^{-1} = g_1$, $\tilde{\sigma} g_3 \tilde{\sigma}^{-1} = g_3g_2$, $f_2 = 2 \epsilon_3 +\epsilon_2$, $\tilde{\sigma} g_4 \tilde{\sigma}^{-1} = g_4^{-1}$, $2\epsilon_3 + 2\delta_4 +\epsilon_2 -1= 0$, $\tilde{\sigma}^2 = Id$, $\gamma_1 =0$ and (\ref{OFG}) splits.\\

1)A1)(b)(i)'' $(E_{\alpha}, \sigma_1)$ as in case A1), $\gamma_2 =0$, $\epsilon_1 =0$. $\tilde{\sigma} g_2 \tilde{\sigma}^{-1} = g_2^{-1}$, $\tilde{\sigma} g_1 \tilde{\sigma}^{-1} = g_1$, $\tilde{\sigma} g_3 \tilde{\sigma}^{-1} = g_3g_2$, $f_2 = 2 \epsilon_3 +\epsilon_2$, $\tilde{\sigma} g_4 \tilde{\sigma}^{-1} = g_4^{-1}$, $2\epsilon_3 + 2\delta_4 +\epsilon_2 -1= 0$, $\tilde{\sigma}^2 = g_1$, $\gamma_1 = \delta_1/2$ and we claim that (\ref{OFG}) does not split.
In fact otherwise there would exist an element $g \in G$ such that $(\tilde{\sigma} g)^2 =1$. But any $g \in G$ can be written as follows: $g = g_4^r g_3^s g_1^l g_2^t$, so we have $(\tilde{\sigma}g)^2 = (\tilde{\sigma} g \tilde{\sigma}^{-1}) \tilde{\sigma} (\tilde{\sigma} g \tilde{\sigma}^{-1}) \tilde{\sigma} = g_4^{-r} g_3^s g_2^s g_1^l g_2^{-t} \tilde{\sigma}(g_4^{-r} g_3^s g_2^s g_1^l g_2^{-t}) \tilde{\sigma}^{-1} \tilde{\sigma}^2= g_4^{-r} g_3^{s} g_4^r g_3^s g_1^{2l+1} g_2^s = g_3^{2s} g_2^{mrs +s} g_1^{2l+1}$. Then $(\tilde{\sigma}g)^2 = 1$ implies $2l +1 =0$, a contradiction. \\

1)A1)(b)(ii)' $(E_{\alpha}, \sigma_1)$ as in case A1), $\gamma_2 =0$, $\epsilon_1 =0$. $\tilde{\sigma} g_2 \tilde{\sigma}^{-1} = g_2^{-1}$, $\tilde{\sigma} g_1 \tilde{\sigma}^{-1} = g_1$, $\tilde{\sigma} g_3 \tilde{\sigma}^{-1} = g_3g_2$, $f_2 = 2 \epsilon_3 +\epsilon_2 $, $\tilde{\sigma} g_4 \tilde{\sigma}^{-1} = g_4^{-1}g_1$, $2\epsilon_3 + 2\delta_4 +\epsilon_2 -1= \delta_1$, $\tilde{\sigma}^2 = Id$, $\gamma_1 =0$ and (\ref{OFG}) splits.\\

1)A1)(b)(ii)'' $(E_{\alpha}, \sigma_1)$ as in case A1), $\gamma_2 =0$, $\epsilon_1 =0$. $\tilde{\sigma} g_2 \tilde{\sigma}^{-1} = g_2^{-1}$, $\tilde{\sigma} g_1 \tilde{\sigma}^{-1} = g_1$, $\tilde{\sigma} g_3 \tilde{\sigma}^{-1} = g_3g_2$, $f_2 = 2 \epsilon_3 +\epsilon_2$, $\tilde{\sigma} g_4 \tilde{\sigma}^{-1} = g_4^{-1}g_1$, $2\epsilon_3 + 2\delta_4 +\epsilon_2 -1= \delta_1$, $\tilde{\sigma}^2 = g_1$, $\gamma_1 =\delta_1/2$ and we claim that(\ref{OFG}) splits. In fact one can easily see that $(\tilde{\sigma} g_4 g_1^{-1})^2 = 1$.\\

2)A1) $(E_{\alpha}, \sigma_1)$ as in case A1), $\gamma_2 =0$, $\epsilon_1 =-1/m$.
$\tilde{\sigma} g_2 \tilde{\sigma}^{-1} = g_2^{-1}$, $\tilde{\sigma} g_1 \tilde{\sigma}^{-1} = g_1 g_2$, $\tilde{\sigma} g_3 \tilde{\sigma}^{-1} = g_3$, $f_2 = 2\epsilon_3$, $\tilde{\sigma} g_4 \tilde{\sigma}^{-1} = g_4^{-1}$, $\epsilon_3 + \delta_4 = 1/2$, $\tilde{\sigma}^2 = Id$, $\gamma_1 =0$ and (\ref{OFG}) splits.\\

1)A2)(a)(i)' $(E_{\alpha}, \sigma_1)$ as in case A2), $\gamma_2 =0$, $\epsilon_1 =0$, $m \equiv 0 \ (mod \ 2)$. $\tilde{\sigma} g_2 \tilde{\sigma}^{-1} = g_2^{-1}$, $\tilde{\sigma} g_1 \tilde{\sigma}^{-1} = g_1$, $\tilde{\sigma} g_3 \tilde{\sigma}^{-1} = g_3$, $f_2 = 2 \epsilon_3$, $\tilde{\sigma} g_4 \tilde{\sigma}^{-1} = g_4^{-1}g_2^{-m/2}$, $\epsilon_3 + \delta_4 = 1/2$, $\tilde{\sigma}^2 = g_3$, $\gamma_1 =(\delta_3 -1/4)/2$ and (\ref{OFG}) does not split.\\

1)A2)(a)(i)'' $(E_{\alpha}, \sigma_1)$ as in case A2), $\gamma_2 =0$,  $\epsilon_1 =0$, $m \equiv 0 \ (mod \ 2)$. $\tilde{\sigma} g_2 \tilde{\sigma}^{-1} = g_2^{-1}$, $\tilde{\sigma} g_1 \tilde{\sigma}^{-1} = g_1$, $\tilde{\sigma} g_3 \tilde{\sigma}^{-1} = g_3$, $f_2 = 2 \epsilon_3$, $\tilde{\sigma} g_4 \tilde{\sigma}^{-1} = g_4^{-1}g_2^{-m/2}$,  $\epsilon_3 + \delta_4 = 1/2$, $\tilde{\sigma}^2 = g_3g_1$, $\gamma_1 =(-1/4 + \delta_1 + \delta_3)/2$ and (\ref{OFG}) does not split.\\

1)A2)(a)(ii)' $(E_{\alpha}, \sigma_1)$ as in case A2), $\gamma_2 =0$, $\epsilon_1 =0$, $m \equiv 0 \ (mod \ 2)$. $\tilde{\sigma} g_2 \tilde{\sigma}^{-1} = g_2^{-1}$, $\tilde{\sigma} g_1 \tilde{\sigma}^{-1} = g_1$, $\tilde{\sigma} g_3 \tilde{\sigma}^{-1} = g_3g_2$, $f_2 = 2 \epsilon_3$, $\tilde{\sigma} g_4 \tilde{\sigma}^{-1} = g_4^{-1}g_1g_2^{-m/2}$, $2\epsilon_3 + 2\delta_4 -1= \delta_1$, $\tilde{\sigma}^2 = g_3$, $\gamma_1 =(\delta_3 -1/4)/2$ and (\ref{OFG}) does not split.\\

1)A2)(a)(ii)'' $(E_{\alpha}, \sigma_1)$ as in case A2), $\gamma_2 =0$, $\epsilon_1 =0$, $m \equiv 0 \ (mod \ 2)$. $\tilde{\sigma} g_2 \tilde{\sigma}^{-1} = g_2^{-1}$, $\tilde{\sigma} g_1 \tilde{\sigma}^{-1} = g_1$, $\tilde{\sigma} g_3 \tilde{\sigma}^{-1} = g_3$, $f_2 = 2 \epsilon_3$, $\tilde{\sigma} g_4 \tilde{\sigma}^{-1} = g_4^{-1}g_1g_2^{-m/2}$, $2\epsilon_3 + 2\delta_4 -1= \delta_1$, $\tilde{\sigma}^2 = g_3g_1$, $\gamma_1 =(\delta_3 -1/4 +\delta_1)/2$ and (\ref{OFG}) does not split.\\

2)A2)(i) $(E_{\alpha}, \sigma_1)$ as in case A2), $\gamma_2 =0$, $\epsilon_1 =-1/m$, $m \equiv 0 \ (mod \ 2)$.
$\tilde{\sigma} g_2 \tilde{\sigma}^{-1} = g_2^{-1}$, $\tilde{\sigma} g_1 \tilde{\sigma}^{-1} = g_1 g_2$, $\tilde{\sigma} g_3 \tilde{\sigma}^{-1} = g_3$, $f_2 = 2\epsilon_3$, $\tilde{\sigma} g_4 \tilde{\sigma}^{-1} = g_4^{-1}g_2^{-m/2}$, $\epsilon_3 + \delta_4 = 1/2$, $\tilde{\sigma}^2 = g_3$, $\gamma_1 =(\delta_3 - 1/4)/2$ and (\ref{OFG}) does not split.\\

2)A2)(ii) $(E_{\alpha}, \sigma_1)$ as in case A2), $\gamma_2 =0$, $\epsilon_1 =-1/m$, $m \equiv 1 \ (mod \ 2)$.
$\tilde{\sigma} g_2 \tilde{\sigma}^{-1} = g_2^{-1}$, $\tilde{\sigma} g_1 \tilde{\sigma}^{-1} = g_1 g_2$, $\tilde{\sigma} g_3 \tilde{\sigma}^{-1} = g_3$, $f_2 = 2\epsilon_3$, $\tilde{\sigma} g_4 \tilde{\sigma}^{-1} = g_4^{-1}g_2^{(1-m)/2}$,  $2\epsilon_3 + 2\delta_4 -1 = \delta_1$, $\tilde{\sigma}^2 = g_3$, $\gamma_1 =(\delta_3 - 1/4)/2$ and (\ref{OFG}) does not split.\\

In order to see that all these cases are different, we only have to check that case 1)A1)(a)(i)' is different from case 1)A1)(a)(i)''; case 1)A1)(b)(ii)' is different from case 1)A1)(b)(ii)''; case 1)A2)(a)(i)' is different from case 1)A2)(a)(i)''; case 1)A2)(a)(ii)' is different from case 1)A2)(a)(ii)''. In fact in each of the above pair of cases, the conjugation action of $\tilde{\sigma}$ is the same, but $\tilde{\sigma}^2$ is different. 

Denote by $\tilde{\sigma}$ and $\tilde{\sigma}'$ the liftings of $\sigma$ as in the above list respectively in the first and in the second case of each pair. Then since the conjugation action of $\tilde{\sigma}$ and of $\tilde{\sigma}'$ on the $g_i's$ are the same, the two cases of each pair are equivalent if and only if there exists a $g \in Z$ such that $\tilde{\sigma}' = \tilde{\sigma} \circ g$. But then if we set $g = g_1^r g_2^t$, $\tilde{\sigma} g \tilde{\sigma}^{-1} = g'$, we have $\tilde{\sigma}'^2 = g' (\tilde{\sigma} g' \tilde{\sigma}^{-1}) \tilde{\sigma}^2$. Now in all these cases we have $\tilde{\sigma}'^2 \tilde{\sigma}^{-2} = g_1$, then we must have $g_1 =    \tilde{\sigma}'^2 \tilde{\sigma}^{-2} = g' (\tilde{\sigma} g' \tilde{\sigma}^{-1}) = g_1^{2r}$, as one can easily compute, so we have found a contradiction.

\hfill{Q.E.D.}

\section{The topology of the real part}

Now we would like to determine the topology of the real part of a real Kodaira surface $S$.
 
First of all we recall that remark (\ref{empty}) tells us that in cases 
1B)'', 1A1)(a)(i)'', 1A1)(b)(i)'', 1)A2)(a)(i)', 1)A2)(a)(i)'',  
1)A2)(a)(ii)', 1)A2)(a)(ii)'', 1)A2)(i), 2)A2)(ii) of \ref{classification}, the real part of $S$ is empty, since (\ref{OFG}) does not split.  

\begin{REM} 
\label{toriklein}
The fixed point locus of an antiholomorphic involution $\sigma$ on a real Kodaira surface $S$ can only be a disjoint union of
tori and Klein bottles.\\  In fact we have the fibration
$$\diagram  S \dto^{\alpha} \rto^{\sigma}    & S \dto^{\alpha}  \\
              E_{\alpha} = {\Co}/({\Z}\alpha_3 + {\Z}/\alpha_4)  \rto^{{\sigma}_1} & E_{\alpha} = {\Co}/({\Z}\alpha_3 + {\Z}/\alpha_4) 
                                        \enddiagram$$ with fibre $E_{\beta} = {\Co}/({\Z}\beta_1 + {\Z}\beta_2)$, so each component of the
real locus is a $S^1$ bundle on $S^1$.

Furthermore the number of such connected components is less or equal to 4, since the real part of an elliptic curve has at most 2 connected components.
\end{REM}

We want now to show how one can compute the number $k$ of connected components of the real part
$S({\R})$ of a real Kodaira surface and how to determine the nature of the components.

Assume that $S({\R}) \neq \emptyset$, then of course $Fix(\sigma_1) \neq \emptyset$ ($\sigma_1$ is the action of $\sigma$ on $E_{\alpha}$ as in \ref{toriklein}) and $Fix(\sigma_1)$ has either one or two connected components homeomorphic to $S^1$.

We consider then all the possible liftings $\tilde{\sigma}'$ of $\sigma$, such that $\tilde{\sigma}'^2 = Id$ and we consider their fixed point loci. Then, for a component in the fixed point locus of such a lifting $\tilde{\sigma}'$, we take the equivalence class, where we say that two such components $\Gamma$ and $\Gamma'$ are equivalent if and only if there exists an element $h \in G$ such that $h (\Gamma) = \Gamma'$. 

Finally we want to see if the components are tori or Klein bottles and we observe that if $\Gamma$ is a connected component of the fixed point locus of a lifting $\tilde{\sigma}'$ as above, the corresponding connected component of $S({\R})$ is homeomorphic to $\Gamma/H$, where $H = \{g \in G \ | g(\Gamma) = \Gamma \}$.

Furthermore, since $\Gamma \cong {\R}^2$, we see that $\pi_1(\Gamma/H) \cong H$, thus we only need to determine $H$.

We will see that $H$ is always abelian, therefore all the connected components are tori.

In order to understand better what we are saying, we show the computation in case 1)A1)(a)(i)' of \ref{classification}.
From theorem \ref{classification} we know that we may choose a lifting $\tilde{\sigma}$ of $\sigma$ to the universal covering ${\R}^4$ such that 
$$\tilde{\sigma} g_2 \tilde{\sigma}^{-1} = g_2^{-1},$$
$$\tilde{\sigma} g_1 \tilde{\sigma}^{-1} = g_1,$$
$$\tilde{\sigma} g_3 \tilde{\sigma}^{-1} = g_3,$$
$$\tilde{\sigma} g_4 \tilde{\sigma}^{-1} = g_4^{-1}.$$

Furthermore we know that $\tilde{\sigma}^2 =1$.

We want to find the elements $g \in G$ such that $(\tilde{\sigma}g)^2 =1$.

We know that every element $g \in G$ can be written as follows: $g = g_4^b g_3^a g_2^t g_1^l$, with $b,a,t,l \in {\Z}$, so we have:
$$(\tilde{\sigma}g)^2 = (\tilde{\sigma} g_4^b g_3^a g_2^t g_1^l \tilde{\sigma}^{-1}) \tilde{\sigma} (\tilde{\sigma} g_4^b g_3^a g_2^t g_1^l \tilde{\sigma}^{-1}) \tilde{\sigma} =$$
$$g_4^{-b} g_3^a g_4^b g_3^a g_1^{2l} = g_4^{-b} g_3^a g_4^b g_3^{-a} g_3^{2a} g_1^{2l} = 
g_2^{mab} g_3^{2a} g_1^{2l}.$$ 

Thus $(\tilde{\sigma}g)^2 =1$ if and only if $a =l=0$, and $g = g_4^b g_2^t$.

Now we want to determine the fixed point locus of $\tilde{\sigma} g$ in the universal covering ${\R}^4$ and we can explicitely compute:
$$\tilde{\sigma} g = \tilde{\sigma}g_4^b g_2^t
\left( \begin{array}{c} 
x_1\\
y_1\\
x_2\\
y_2 
\end{array} \right)
= \left( \begin{array}{cccc} 
1 & 0 & 0 & 0\\
0 & -1 & 0 & 0\\
0 & f_2 +b & 1 & 0\\
f_2+b & 0 & 0 & -1 
\end{array} \right)
\left( \begin{array}{c} 
x_1\\
y_1\\
x_2\\
y_2 
\end{array} \right)+ 
\left( \begin{array}{c} 
0\\
-b\\
f_2 b + b(b-1)/2 + b \delta_4\\
-b \epsilon_4 -t \epsilon_2 
\end{array} \right).
$$
So $\tilde{\sigma} g (x) = x$ if and only if 
$$y_1 = -b/2,$$
$$(f_2 + b)y_1 + f_2 b + b(b-1)/2 + b \delta_4 =0,$$
$$y_2 = (f_2 + b)x_1/2 - b\epsilon_4/2 -t\epsilon_2/2.$$
Now one can easily prove that we can assume $b=0,1$ and we find at most 4 connected components.
In fact we can see that if we take $b, b'$ as above such that $b \equiv b' \ (mod \ 2)$, and we call $\Gamma$, $\Gamma'$ the corresponding components, there exists an element $h \in G$ such that $h(\Gamma) = \Gamma'$.

Therefore it suffices to make the explicit computation for $b =0,1$.

If $b =0$ we have $y_1 =0$ and the third equation above gives $y_2 = \frac{f_2}{2} x_1  - t \frac{\epsilon_2}{2}$ and we get two connected components 
$$\Gamma_1 = \{y_1 =0, \ y_2 = \frac{f_2}{2} x_1\},$$
$$\Gamma_2 = \{y_1 =0, \ y_2 = \frac{f_2}{2} x_1 + \frac{\epsilon_2}{2} \}.$$
Now we must check whether the two components give rise to different components of $S({\R})$ or if there exists an element $h \in G$ such that $h(\Gamma_1) = \Gamma_2$. 

If $b = 1$, from the second equation above we obtain $f_2/2 - 1/2  + \delta_4 = 0$ and this holds trivially since we have $\delta_4 = 1/2 - f_2/2$ (cf. \ref{classification}).

So we finally get at most two other connected components, namely 
$$\Lambda_1 = \{y_1 = -1/2, \ y_2 =(\frac{f_2 +1}{2}) x_1 - \frac{\epsilon_4}{2} \},$$
$$\Lambda_2 = \{y_1 = -1/2, \ y_2 =(\frac{f_2 +1}{2}) x_1 - \frac{\epsilon_4}{2} + \frac{\epsilon_2}{2} \}.$$
Again one has to check whether these two components give rise to two different connected components in $S({\R})$ or not.

We show that $\Gamma_1$ and $\Gamma_2$ are not in the same equivalence class.
In fact assume that there exists a $g \in G$ such that $g(\Gamma_1) = \Gamma_2$, then we can write $g = g_4^s g_3^r g_1^l g_2^n$ and we would have:

$$
\left( \begin{array}{cccc} 
1 & 0 & 0 & 0\\
0 & 1 & 0 & 0\\
r & s & 1 & 0\\
-s & r & 0 & 1 
\end{array} \right)
\left( \begin{array}{c} 
x_1\\
0\\
x_2\\
\epsilon_3 x_1 
\end{array} \right)+ 
\left( \begin{array}{c} 
r\\
s\\
\delta\\
\epsilon 
\end{array} \right)=
\left( \begin{array}{c} 
x'_1\\
0\\
x'_2\\
\epsilon_3 x'_1 + \epsilon_2/2 
\end{array} \right),$$
so $s =0$ and $g = g_3^r g_1^l g_2^n$, which yields $\epsilon = r \epsilon_3 + n \epsilon_2$ (recall that by \ref{classification} we have $f_2 = 2 \epsilon_3$). Then we must have $x'_1 = x_1 + r$, $\epsilon_3 x_1 + \epsilon = \epsilon_3 x'_1 + \epsilon_2/2 =  \epsilon_3 (x_1 + r) + \epsilon_2/2$, so we get $n \epsilon_2 = \epsilon_2/2$, which implies $n = 1/2$, a contradiction. 

Therefore there cannot exist a $g \in G$ such that $g (\Gamma_1) = \Gamma_2$, so $\Gamma_1$ and $\Gamma_2$ are not equivalent and they give rise to two distinct components of $S({\R})$.

For the $\Lambda_i$'s, we show that if $m \equiv 0 \ (mod \ 2)$, then $\Lambda_1$ and $\Lambda_2$ are not equivalent, while if $m \equiv 1 \ (mod \ 2)$, they are equivalent.

In fact we want to see if there exists an element $g \in G$ such that $g (\Lambda_1) = \Lambda_2$, so we write $g = g_4^s g_3^r g_1^l g_2^n$ as above and we impose

$$
\left( \begin{array}{cccc} 
1 & 0 & 0 & 0\\
0 & 1 & 0 & 0\\
r & s & 1 & 0\\
-s & r & 0 & 1 
\end{array} \right)
\left( \begin{array}{c} 
x_1\\
-1/2\\
x_2\\
(\epsilon_3 + 1/2)x_1 - \epsilon_4/2 
\end{array} \right)+ 
\left( \begin{array}{c} 
r\\
s\\
\delta\\
\epsilon 
\end{array} \right)=
\left( \begin{array}{c} 
x'_1\\
-1/2\\
x'_2\\
(\epsilon_3 + 1/2)x'_1 -\epsilon_4/2 + \epsilon_2/2 
\end{array} \right),$$
so $s =0$ and $g = g_3^r g_1^l g_2^n$, which yields $\epsilon = r \epsilon_3 + n \epsilon_2$. Then we must have $x'_1 = x_1 + r$, $-r/2 +(\epsilon_3 +1/2) x_1 -\epsilon_4/2 + \epsilon = (\epsilon_3 +1/2) x'_1 - \epsilon_4/2 + \epsilon_2/2 =  (\epsilon_3 + 1/2)(x_1 + r) -\epsilon_4/2 + \epsilon_2/2$, so we get $n \epsilon_2 - \epsilon_2/2 = r$. Since $\epsilon_2 = 2/m$, we obtain $2n -1 = mr$, which is obviously impossible if $m \equiv 0 \ (mod \ 2)$. If instead $m \equiv 1 \ (mod \ 2)$, we can for instance choose $r =1$, $n = (m+1)/2$, so $\Lambda_1$ and $\Lambda_2$ are equivalent. 

Therefore we have proven that if $m \equiv 0 \ (mod \ 2)$, $S({\R})$ has 4 distinct connected components, while if $m \equiv 1 \ (mod \ 2)$, $S({\R})$ has 3 distinct connected components.

Now we show that the component $\Gamma_1$ gives rise to a torus as a component of $S({\R})$.

We have then to determine the group $H = \{g \in G \ | \ g(\Gamma_1) = \Gamma_1 \}$. Since any $g \in G$ is of the form

$$g 
\left( \begin{array}{c} 
x_1\\
y_1\\
x_2\\
y_2 
\end{array} \right)
= \left( \begin{array}{cccc} 
1 & 0 & 0 & 0\\
0 & 1 & 0 & 0\\
r & s & 1 & 0\\
-s & r & 0 & 1 
\end{array} \right)
\left( \begin{array}{c} 
x_1\\
y_1\\
x_2\\
y_2 
\end{array} \right)+ 
\left( \begin{array}{c} 
r\\
s\\
\delta\\
\epsilon 
\end{array} \right),
$$
and $\Gamma_1 = \{y_1 =0, \ y_2 = \frac{f_2}{2} x_1 \}$, where $f_2 = 2 \epsilon_3$ (cf. \ref{classification}), we want 
$$
\left( \begin{array}{cccc} 
1 & 0 & 0 & 0\\
0 & 1 & 0 & 0\\
r & s & 1 & 0\\
-s & r & 0 & 1 
\end{array} \right)
\left( \begin{array}{c} 
x_1\\
0\\
x_2\\
\epsilon_3 x_1 
\end{array} \right)+ 
\left( \begin{array}{c} 
r\\
s\\
\delta\\
\epsilon 
\end{array} \right)=
\left( \begin{array}{c} 
x'_1\\
0\\
x'_2\\
\epsilon_3 x'_1 
\end{array} \right).$$
So we find $x'_1 = x_1 + r$, $s =0$, $g= g_3^r g_1^n g_2^t$ ($ r,n,t \in {\Z}$), $\epsilon_3 x_1 + \epsilon = \epsilon_3 (x_1 + r)$. Now since $g = g_3^r g_1^n g_2^t$, we obtain $\epsilon = \epsilon_3 r + t \epsilon_2$, so we must have $t =0$, $g = g_3^r g_1^n$, $r, n \in {\Z}$. So $H = \{g_3^rg_1^n, \ | \ r, n \in {\Z}\} \cong {\Z}^2$ and $\Gamma_1/H \cong S^1 \times S^1$.

Analogously we can prove that also the other components are tori.

Similar computations allow us to describe the topology of the real part in all the other cases and we obtain the following table.

$$\begin{tabular}{|c|c|c|}
\hline
$Case$ & $m$ & $S({\R})$  \\
\hline
\hline
$1)B)'$ &  $m \equiv 0 \ (mod \ 2)$ & $2T$ \\

\hline
$1)B)'$ & $m \equiv 1 \ (mod \ 2)$  & $T$\\
\hline

$2)B)$ & $\forall m $ & $T$\\

\hline

$1)A1)(a)(i)'$ & $m \equiv 0 \ (mod \ 2)$ &  $4T$\\

\hline
$1)A1)(a)(i)'$ & $m \equiv 1 \ (mod \ 2)$ &  $3T$\\

\hline

$1)A1)(a)(ii)'$ & $ \forall m $ &  $2T$\\
\hline

$1)A1)(a)(ii)''$ & $m \equiv 0 \ (mod \ 2)$ &  $2T$\\

\hline
$1)A1)(a)(ii)''$ & $m \equiv 1 \ (mod \ 2)$ &  $T$\\

\hline
$1)A1)(b)(i)'$ & $m \equiv 0 \ (mod \ 2)$ &  $3T$\\

\hline
$1)A1)(b)(i)'$ & $m \equiv 1 \ (mod \ 2)$ &  $4T$\\

\hline

$1)A1)(b)(ii)'$ & $\forall m$ &  $T$\\

\hline

$1)A1)(b)(ii)''$ & $m \equiv 0 \ (mod \ 2)$ &  $T$\\

\hline
$1)A1)(b)(ii)''$ & $m \equiv 1 \ (mod \ 2)$ &  $2T$\\

\hline
$2)A1)$ & $ \forall m $ &  $2T$\\
\hline 
$1)B)''$ & $\forall m$ &  $\emptyset$\\

\hline 
$1)A1)(a)(i)''$ & $\forall m$ &  $\emptyset$\\
\hline

$1)A1)(b)(i)''$ & $\forall m$ &  $\emptyset$\\
\hline

$1)A2)(a)(i)'$ & $m \equiv 0 \ (mod \ 2)$ &  $\emptyset$\\
\hline

$1)A2)(a)(i)''$ & $m \equiv 0 \ (mod \ 2)$ &  $\emptyset$\\
\hline

$1)A2)(a)(ii)'$ & $m \equiv 0 \ (mod \ 2)$ &  $\emptyset$\\
\hline
$1)A2)(a)(ii)''$ & $m \equiv 0 \ (mod \ 2)$ &  $\emptyset$\\
\hline

$2)A2)(i)$ & $m \equiv 0 \ (mod \ 2)$ &  $\emptyset$\\
\hline
$2)A2)(ii)$ & $m \equiv 1 \ (mod \ 2)$ &  $\emptyset$\\
\hline 
\end{tabular}
$$

\bigskip

Address of the author:\\

Paola Frediani: \\ Dipartimento di Matematica della Universit\`a di Pavia\\
via Ferrata, 1, I-27100 Pavia

e-mail: frediani@dimat.unipv.it\\


\begin{thebibliography} {Bellll,Fr,Sol}

\bibitem[A-G]{ag} N. L. Alling, N. Greenleaf ``{\em Foundations of the
theory of Klein surfaces}'', Volume 219 of {\em Lecture Notes in
Mathematics}, Springer-Verlag, Berlin-Heidelberg-New York, (1980).


\bibitem[B-P-VdV]{bpvdv} W. Barth, C. Peters, A. Van de Ven, ``{\em Compact
Complex Surfaces}'', Springer Ergebnisse 3, 4, Springer - Verlag,(1984).

\bibitem[Be]{be} A. Beauville ``{\em  Surfaces alg\'ebriques complexes}'',
Ast\'erisque, No. 54. Soci\'et\'e Math\'ematique de France, Paris (1978).

\bibitem[Bo1]{bo1} C. Borcea ''{\em Moduli for Kodaira surfaces},
Comp. Math. {\bf 52} 373-380 (1984).


\bibitem[Cal]{cal} E. Calabi ''{\em Construction and properties of some
6-dimensional almost complex manifolds}, Trans. Amer.  Math. Soc. {\bf
187} 407-438 (1958).

\bibitem[Ca1]{ca1} F. Catanese ''{\em On the moduli spaces of surfaces of
general type},  Jour. Diff. Geom.  {\bf 19} 483-515 (1984).

\bibitem[Ca2]{ca2} F. Catanese ''{\em Connected components of moduli spaces},
Jour. Diff. Geom.  {\bf 24} 395-399 (1986).


\bibitem[Ca5]{ca5} F. Catanese ''{\em Moduli spaces of surfaces and real structures}, 
AG/0103071 (2001), to appear in Annals of Math.

\bibitem[Ca]{ca} F. Catanese ''{\em Fibred surfaces, varieties isogenous to a
product and related moduli spaces}, American Journal of Math. {\bf 122} 1-44
(2000).

\bibitem[C-F]{cf} F. Catanese, P. Frediani ''{\em Real hyperelliptic surfaces and the orbifold fundamental group},  to appear in Journal de l'Institute de Math. de Jussieu, 2 (2003).  

\bibitem[Ci-Pe]{ci-pe} C.Ciliberto, C. Pedrini ``{\em Real abelian varieties
and real algebraic curves}'' Lectures in real geometry (Madrid, 1994),
Fabrizio Broglia ed., de Gruyter Exp. Math., 23, Berlin (1996) 167--256.

\bibitem[CO1]{co1} A. Comessatti ``{\em Fondamenti per la geometria sopra le
superficie razionali dal punto di vista reale}'', Math. Annalen {\bf 43} 1-72
(1912).

\bibitem[CO2]{co2} A. Comessatti ``{\em Sulla connessione delle superficie
razionali reali}'', Annali di Matemetica {\bf 23} 215-285 (1914).

\bibitem[CO3]{co3} A. Comessatti ``{\em Reelle Fragen in der algebraischen
Geometrie}'', Jahresbericht d. Deut. Math. Vereinigung {\bf 41} 107-134
(1932).

\bibitem[CO4]{co4} A. Comessatti ``{\em Sulle variet\`a abeliane reali I e II}, Ann. Mat. Pura Appl. 2 and 4, 67-106 and 27-71,  (1924 and 1926).

\bibitem[D-I-K]{dik} A. Degtyarev, I. Itenberg, V. Kharlamov ''{\em Real Enriques
Surfaces} Springer Lecture Notes 1746, (2000).

\bibitem[D-K1]{dekha1} A. Degtyarev, V. Kharlamov ''{\em Topological
classification of real Enriques surfaces}'', Topology {\bf 35} 711-729 (1996).

\bibitem[D-K2]{dekha2} A. Degtyarev, V. Kharlamov ''{\em Topological
properties of real algebraic varieties: de cot\'e de chez Rokhlin}'', Russ. Math. Surveys. 55, no.4, 735-814 (2000).

\bibitem[D-K3]{dekha3} A. Degtyarev, V. Kharlamov ''{\em On the moduli space
of real Enriques surfaces}'', C. R. Acad. Sci. Paris S\'er. I Math. 324, 3,
317-322 (1997).

\bibitem[D-K4]{dekha4} A. Degtyarev, V. Kharlamov ''{\em Real rational surfaces are quasi-simple}'', J. Reine Angew. Math. 551, 87-99, (2002). 




\bibitem [Kha]{kha} V. Kharlamov ''{\em The topological type of nonsingular
surfaces in ${\R}{\proj}^3 $ of degree $4$}'', Functional Anal. Appl.{\bf 10}
295-305 (1976).

\bibitem [Kha-Ku]{kha-ku} V. Kharlamov, Vik. S. Kulikov  ''{\em On 
real structures on rigid surfaces}'', Izv. Ross. Akad. Nauk Ser. Mat. 66 (2002), no. 1, 133-152; translation in Izv. Math. 66 (2002), no. 1, 133-150.

\bibitem [K-M]{k-m} K. Kodaira, J. Morrow ''{\em Complex manifolds} Holt,
Rinehart and Winston, Inc., New York-Montreal, Que.-London (1971).


\bibitem [Ko]{ko} K. Kodaira ''{\em Complex manifolds and deformation of
complex structures} Springer  Grundlehren der Math. Wiss.
{\bf 283},   New
York-Berlin, (1986).

\bibitem [KoI]{koi} K. Kodaira ''{\em On the structure of complex analytic
surfaces I}'', Amer. J. Math. {\bf 90} 1048--1066 (1968).

\bibitem [KoIV]{koiv} K. Kodaira ''{\em On the structure of complex analytic
surfaces IV}'', Amer. J. Math. {\bf 90} 1048--1066 (1968).

\bibitem[Ma]{ma} M. Manetti ''{\em On the moduli space of diffeormorphic
algebraic surfaces } Inv. math.,  to appear

\bibitem [Man1]{man1} F. Mangolte''{\em Cycles alg\'ebriques sur les surfaces
$K3$ r\'eelles}''.  Math. Z. 225, 4, 559--576 (1997).

\bibitem [Man2]{man2} F. Mangolte, J. van Hamel''{\em Algebraic cycles and
topology of real Enriques surfaces}''. Compositio Math. 110, 2,
215--237(1998).

\bibitem [Man3]{man3} F. Mangolte''{\em Surfaces elliptiques r\'eelles et
in\'egalit\'es de Ragsdale-Viro}'',  Math. Z. 235 (2000), no. 2, 213--226. 

\bibitem [Ni]{ni} V.V. Nikulin ''{\em Integral symmetric bilinear forms and
some of their applications}'', Izv. Akad. Nauk SSSR {\bf 43} 1, 117-177
(1979) (Russian); English transl. in Math. USSR-Izv. {\bf 14} 103-167 (1980).

\bibitem [Se]{mika} M. Sepp\"al\"a ''{\em Real algebraic curves in the moduli
space of complex curves}''. Compositio Math. 74,3, 259--283 (1990).

\bibitem [Se-Si]{s-s} M. Sepp\"al\"a, R. Silhol ''{\em Moduli spaces for real
algebraic curves and real abelian varieties}''. Math. Z. 201, 2, 151--165
(1989).


\bibitem[Si]{si} R. Silhol ``{\em Real Algebraic Surfaces}'', Lectures Notes
in Mathematics 1392, Springer - Verlag (1989).

\bibitem[We]{we} Jean-Yves Welschinger ``{\em Real structures on minimal ruled surfaces}'', Comment. Math. Helv. 78, 418-446 (2003).




\end{thebibliography}
\end{document}